\magnification=1000
\hsize=11.7cm
\vsize=18.9cm
\lineskip2pt \lineskiplimit2pt
\nopagenumbers

\hoffset=-1truein
\voffset=-1truein

\advance\voffset by 4truecm
\advance\hoffset by 4.5truecm

\newif\ifentete

\headline{\ifentete\ifodd	\count0 
      \rlap{\head}\hfill\tenrm\llap{\the\count0}\relax
    \else
        \tenrm\rlap{\the\count0}\hfill\llap{\head} \relax
    \fi\else
\global\entetetrue\fi}

\def\entete#1{\entetefalse\gdef\head{#1}} 
\entete{}

\input amssym.def
\input amssym.tex

\def\-{\hbox{-}}
\def\.{{\cdot}}

\def\F{{\cal F}}

\def\L{{\cal L}}
\def\M{{\cal M}}

\def\P{{\cal P}}

\def\T{{\cal T}}

\def\C{{\cal C}}

\def\Ab{\frak A\frak b}

\def\Gr{\frak G\frak r}

\def\Fct{\frak F\frak c\frak t}

\def\Ker{\frak K\frak e\frak r}

\def\int{\frak i\frak n\frak t}

\def\qq{\quad{\rm and}\quad}

\def\too{\longrightarrow}

 3
 2
\font\large=cmr10  scaled \magstep 2
 2
\font\larti=cmti10  scaled \magstep 2
 1
 2

\font\cds=cmr7
\font\cdt=cmti7
\font\cdy=cmsy7

\count0=1

\centerline{\large A correction to the uniqueness of a partial}
\medskip
\centerline{\large  perfect locality over a Frobenius {\larti P}-category}
\medskip
\centerline{\bf Lluis Puig }
\medskip
\noindent 
\centerline{\cds CNRS, Institut de Math\'ematiques de Jussieu, lluis.puig@imj-prg.fr}
\par
\noindent
\centerline{\cds 6 Av Bizet, 94340 Joinville-le-Pont, France}

\medskip
\noindent
{\bf Abstract:} {\cds Let {\cdt p} be a prime, {\cdt P} a finite {\cdt p}-group and {\cdy F} a {\cdt Frobenius P-category\/}.  In  {\cdt Existence, uniqueness and functoriality of the perfect locality over a Frobenius 
P-category}, Algebra Colloquium, 23(2016) 541-622, we also claimed the uniqueness of the partial 
perfect locality {\cdy L}{${}^{_{\frak X}}$} over  any up-closed set $\scriptstyle \frak X$ of {\cdt {\cdy F}-selfcentralizing\/} subgroups of {\cdt P}, but recently Bob Oliver exhibit some counter-examples, demanding
some revision of our arguments. In this Note we show that, up to replacing the perfect localities by the {\cdt extendable\/} perfect localities over  any up-closed set $\scriptstyle \frak X$ of {\cdt {\cdy F}-selfcentralizing\/} subgroups of {\cdt P}, our arguments are correct, still proving the existence and the uniqueness of the {\cdt perfect {\cdy F}{${}^{_{\rm sc}}\!$}-locality\/}, since it is {\cdt extendable\/}.
We take advantage to simplify some of our arguments.}

\bigskip
\noindent
{\bf £1. Introduction  }
\medskip
£1.1.  Let $p$ be a prime, $P$ a finite $p\-$group, $\F$ a {\it Frobenius  $P\-$category\/}~[2]
and $\T_P$ the category where the objects are the subgroups of $P\,,$ the morphisms are defined by the
{\it $P\-$transporters\/} and the composition is defined by the product in $P\,.$
Recall that, according to~[3,~17.3], an {\it $\F\-$locality\/} $\L$ is a finite category where the objects are all  the subgroups of~$P\,,$ endowed with two functors
$$\tau : \T_{\! P}\too \L\qq \pi : \L\too \F
\eqno £1.1.1\phantom{.}$$
which are the identity on the set of objects, $\pi$ being {\it full\/}, and such that the composition 
$\pi\circ\tau$ is induced by the conjugation in $P\,;$ we say that $\L$ is {\it divisible\/} whenever it fulfills the following condition
\smallskip
\noindent
£1.1.2\quad {\it If $Q\,,$ $R$ and $T$ are subgroups of $P\,,$ for any $\L\-$morphisms
$x\,\colon R\to Q$ and $y\,\colon T\to Q$ such that the image of $\pi_{_{Q,T}}(y)$ is contained in
 the image of $\pi_{_{Q,R}}(x)\,,$ there is a
unique $\L\-$morphism $z\,\colon T\to R$ such that $x\.z = y\,.$\/}

\medskip
£1.2. Then, it follows from [3,~Proposition~18.4 and~Theorem~18.6] that a 
{\it perfect $\F\-$locality\/}, introduced in~[3,~17.13], is a {\it divisible $\F\-$locality\/} $\P$ such that, for any subgroup $Q$ of $P$ {\it fully normalized\/} in $\F$ [3,~2.6] the finite group $\P (Q)$ endowed with the group homomorphims
$$\tau_{_Q} : N_P(Q)\too \P (Q)\qq \pi_{_Q} : \P (Q)\too \F (Q)
\eqno £1.2.1$$
is the {\it $\F\-$localizer\/} of $Q\,,$ introduced in~[3,~18.5]. Actually, as we show in [3,~Theorem~20.24] and, more carefully, in [5,~Theorem~7.2], $\P$ is  uniquely determined by
the {\it full\/} subcategory $\P^{^{\rm sc}}$ over the set of {\it $\F\-$selfcentralizing\/} subgroups
of $P\,,$ introduced in~[3,~4.8].

\medskip
£1.3. More generally, in order to apply inductive arguments, for any nonempty set $\frak X$ of 
$\F\-$selfcentralizing  subgroups of $P$ which contains any\break
\eject
\noindent 
subgroup of $P$ admitting an $\F\-$morphism from some subgroup in~$\frak X\,,$ we consider the {\it full\/} subcategory 
$\F^{^\frak X}$ of $\F$ over $\frak X$ as the set of objects and, replacing $\T_P$ by its {\it full\/}
subcategory $\T^{^\frak X}_P$ over $\frak X\,,$ we may introduce the {\it $\F^{^\frak X}\-$localities\/}
as the finite categories $\L^{^\frak X}$ where the objects are  the groups in~$\frak X\,,$ endowed with two functors
$$\tau^{_\frak X} : \T^{^\frak X}_{\! P}\too \L^{^\frak X}\qq \pi^{_\frak X} : \L^{^\frak X}\too \F^{^\frak X}
\eqno £1.3.1\phantom{.}$$
which are the identity on the set of objects, $\pi^{_\frak X}$ being {\it full\/}, and such that the composition $\pi^{_\frak X}\circ\tau^{_\frak X}$ is induced by the conjugation in $P\,.$
In particular, in~[5,~2.8] we consider a {\it perfect $\F^{^\frak X}\!\-$locality\/} $\P^{^\frak X}$
and in~[5,~6.1] we claimed its existence and uniqueness.

\medskip
£1.4.  But, recently, Bob Oliver exhibit some counter-examples to this uniqueness [1]; of course, these counter-examples demand a revision of our arguments
in~[5]. Our purpose in this Note is to show that, up to restricting the {\it perfect 
$\F^{^\frak X}\!\-$localities\/} we consider, our arguments become correct and the 
uniqueness of these restricted {\it perfect $\F^{^\frak X}\!\-$localities\/}, called {\it extendable\/},  is true; naturally, our {\it extendable perfect $\F^{^{\rm sc}}\!\-$localities\/}
include $\P^{^{\rm sc}}$ above. Moreover, we take advantage of this revision to simplify some arguments in~[5]. Notations and terminology are the same as in~[5] and the main 
references come from~[3].

\bigskip
\noindent
{\bf £2. Extendable perfect $\F^{^\frak X}\-$localities\/}

\medskip
£2.1. With the notation above, let us consider a {\it perfect $\F^{^\frak X}\!\-$locality\/} 
$\P^{^\frak X}\,;$ that is to say, $\P^{^\frak X}$ is a {\it divisible $\F^{^\frak X}\!\-$locality\/}
 such that, for any group $Q$ in~$\frak X$ {\it fully normalized\/} in $\F$ [3,~2.6], the finite group 
 $\P^{^\frak X}\! (Q)$ endowed with the group homomorphims
$$\tau^{_\frak X}_{_Q} : N_P(Q)\too \P^{^\frak X} \!(Q)\qq \pi^{_\frak X}_{_Q} : \P^{^\frak X}\! (Q)\too \F (Q)
\eqno £2.1.1$$
is the {\it $\F\-$localizer\/} of $Q\,,$ introduced in~[3,~18.5]; in particular, $\pi^{_\frak X}_{_Q}$
is surjective and, since $Q$ is $\F\-$selfcentralizing, $\tau^{_\frak X}_{_Q}$ is injective [3,~Remark~18.7]. Moreover, note that  condition~18.6.3 in [3,~Theorem~18.6] implies that
$Q$ fulfills equality~17.10.1 in [3,~Proposition~17.10]; in particular, extending $\P^{^\frak X}$
as in [3,~17.4], it follows from [3,~Proposition~17.10] that $\P^{^\frak X}$ is a {\it coherent
$\F^{^\frak X}\!\-$locality\/} [3,~17.9]; that is to say, we have
$$ x\.\tau_{_R} (v) = \tau_{_Q}\big(\pi_x (v)\big)\.x 
\eqno £2.1.2\phantom{.}$$ 
for any pair of subgroups $Q$ and $R$ in $\frak X\,,$ any $x\in\P^{^\frak X} (Q,R)$ and any 
$v\in
R\,.$

\medskip
£2.2. Actually, considering the 
{\it normalizer\/} $N_\F (Q)$ of $Q$ in $\F$ [3,~2.14]  which is a Frobenius $N_P (Q)\-$category
[3,~Proposition~2.16], denoting by $\frak X^{^Q}$ the set of subgroups of $N_P (Q)$ belonging to 
$\frak X$ and setting $P^{^Q} = N_P (Q)$\break
\eject
\noindent
 and $\F^{^\frak Q}\! = N_\F (Q)\,,$  we can also consider
the {\it normalizer\/} $N_{\P^{^\frak X}} (Q)$ of $Q$ in $\P^{^\frak X}$ [3,~17.4 and~17.5]
and, setting $\F^{^{\frak X,Q}} =(\F^{^Q})^{^{\frak X^{^Q}}}\,,$  it is not difficult to see that 
$\P^{^{\frak X,Q}}\!\! = N_{\P^{^\frak X}} (Q)$ is actually a
{\it perfect $\F^{^{\frak X,Q}}\!\-$locality\/}.

\medskip
£2.3. Moreover, the $\F^{^{\frak X,Q}}\!\-$locality $\P^{^{\frak X,Q}}\!$ and the group 
$\P^{^\frak X} \!(Q)$ are related throughout the {\it transporter\/} of the $p\-$subgroups of 
$\P^{^\frak X} \!(Q)\,;$ explicitly, let us call 
{\it transporter\/} $\T_{\P^{^\frak X} \!(Q)}$   of $\P^{^\frak X} \!(Q)$ the $\F^{^Q}\!\-$locality formed by the category where the~objects are all the subgroups of $P^{^Q}\,,$ where the morphisms are defined~by the 
{\it elements\/} of the  $\P^{^\frak X} \!(Q)\-$transporters of the corresponding 
$\tau^{_\frak X}_{_Q}\-$images, and where the composition is defined by the product in $\P^{^\frak X} \!(Q)\,,$ endowed with the obvious functors induced by $\tau^{_\frak X}_{_Q}$ and by~$\pi^{_\frak X}_{_Q}\,.$ Then, denoting by $\T_{\P^{^\frak X} \!(Q)}^{^{\frak X^{^Q}}}$ the {\it full\/} subcategory of 
$\T_{\P^{^\frak X} \!(Q)}$ over $\frak X^{^Q}\,,$  we claim that we have an  
{\it $\F^{^Q}\-$locality equivalence\/} [3,~2.9]
$$\P^{^{\frak X,Q}}\cong \T_{\P^{^\frak X} \!(Q)}^{^{\frak X^{^Q}}}
\eqno  £2.3.1.$$
Firstly, we need the following lemma which admits the same proof as in [3,~Proposition~24.2].

\bigskip
\noindent
{\bf Lemma~£2.4.} {\it Any $\P^{^\frak X}\-$morphism is a monomorphism and an epimorphism.}

\medskip
£2.5. Now, we already know that any $\P^{^{\frak X,Q}}\!\-$morphism $x\,\colon T\to R$ is induced by a
$\P^{^{\frak X}}\!\-$morphism $\hat x\,\colon T\.Q\to R\.Q$ which stabilizes $Q$ [3,~2.14.1];
then, it easily follows from the lemma above that $\hat x$ is uniquely determined by $x\,,$ and the 
{\it divisibility\/} of~$\P^{^{\frak X}}\!$ guarantees the existence of a unique 
$\hat x_{_Q}\in \P^{^{\frak X}}\! (Q)$ fulfilling
$$\tau^{_\frak X}_{_{R\.Q,Q}} (1) \. \hat x_{_Q} = \hat x\. \tau^{_\frak X}_{_{T\.Q,Q}} (1) 
\eqno £2.5.1;$$
moreover, from the {\it coherence\/} of $\P^{^{\frak X}}\!$ (cf.~£2.1.2), for any $t\in T\i P^{^Q}$ we get
$$\eqalign{\tau^{_\frak X}_{_{R\.Q,Q}} (1) \. \hat x_{_Q}\. \tau^{_\frak X}_{_Q} (t) &= 
\hat x\. \tau^{_\frak X}_{_{T\.Q,Q}} (t) = \hat x\. \tau^{_\frak X}_{_{T\.Q}} (t)\. 
\tau^{_\frak X}_{_{T\.Q,Q}} (1)\cr
& = \tau^{_\frak X}_{_{R\.Q}} \Big(\big(\pi^{_\frak X}_{_{R\.Q,T\. Q}}(\hat x)\big)(t)\Big)
\.\hat x\. \tau^{_\frak X}_{_{T\.Q,Q}} (1)\cr
&= \tau^{_\frak X}_{_{R\.Q}} \Big(\big(\pi^{_\frak X}_{_{R\.Q,T\. Q}}(\hat x)\big)(t)\Big)
\.\tau^{_\frak X}_{_{R\.Q,Q}} (1) \. \hat x_{_Q}\cr
&= \tau^{_\frak X}_{_{R\.Q,Q}} (1) \. \tau^{_\frak X}_{_Q} \Big(\big(\pi^{_\frak X}_{_{R\.Q,T\. Q}}(\hat x)\big)(t)\Big)\. \hat x_{_Q}\cr
&= \tau^{_\frak X}_{_{R\.Q,Q}} (1) \. \tau^{_\frak X}_{_Q} \Big(\big(\pi^{_\frak X}_{_{R,T}}( x)\big)(t)\Big)\. \hat x_{_Q}\cr}
\eqno £2.5.2,$$
so that from the lemma above we still get
$$\hat x_{_Q}\. \tau^{_\frak X}_{_Q} (t) \. (\hat x_{_Q})^{-1} = \tau^{_\frak X}_{_Q} \Big(\big(\pi^{_\frak X}_{_{R,T}}( x)\big)(t)\Big)
\eqno £2.5.3;$$
thus, the element $\hat x_{_Q}$ belongs to the $\P^{^{\frak X}}\! (Q)\-$transporter 
$\T_{\P^{^{\frak X}}\! (Q)}\big(\tau^{_\frak X}_{_Q} (R),\tau^{_\frak X}_{_Q} (T)\big)$\break
\eject 
\noindent
and it is not difficult to check that the correspondence sending the $\P^{^{\frak X,Q}}\!\-$mor-phism 
$x\,\colon T\to R$ to the $\T_{\P^{^\frak X} \!(Q)}\-$morphism $\hat x_{_Q}\,\colon T\to R$
defines a {\it faithful $\F^{^Q}\-$locality functor\/} $\P^{^{\frak X,Q}}\to 
\T_{\P^{^\frak X} \!(Q)}^{^{\frak X^{^Q}}}$ [3,~2.9]. The ``surjectivity'' follows again from
condition~18.6.3 in [3,~Theorem~18.6].

\medskip
£2.6. But, for any $\F\-$selfcentralizing subgroup $W$ of $P$ {\it fully normalized\/} in~$\F\,,$
we still have the normalizer $\F^{^W} = N_\F (W)\,;$ let us set $P^{^W} = N_P (W)\,;$ 
if $P^{^W}$ belongs to $\frak X\,,$ so that  the set $\frak X^{^W}$ of subgroups of $P^{^W}$ belonging to~$\frak X$ is not empty, then we also
can consider the normalizer $\P^{^{\frak X,W}} \!= N_{\P^{^\frak X}} (W)\,,$ which is again a
{\it perfect $\F^{^{\frak X,W}}\!\-$locality\/}, and we always have the existence of the
{\it $\F\-$localizer\/} $L^{^W}_\F$ of $W$ [3,~Theorem~18.6]; thus, we still can consider the 
{\it transporter\/} $\T_{L^{^W}_\F}$ of~$L^{^W}_\F$ as an $\F^{^W}\!\-$locality and  the {\it full\/} subcategory  $\T_{L^{^W}_\F}^{^{\frak X^{^W}}}\!$ of~$\T_{L^{^W}_\F}$
over $\frak X^{^W}\!$ as the set of objects. Finally, we say that the  {\it perfect $\F^{^\frak X}\!\-$locality\/} $\P^{^\frak X}$ is {\it extendable\/} whenever for any $\F\-$selfcentralizing subgroup $W$ of $P$ {\it fully normalized\/} in~$\F$ such that $P^{^W}\in \frak X$ there exists an {\it $\F^{^{\frak X,W}}\!\-$locality isomorphism\/}{\footnote{\dag}{
\cds In [5,~6.18], arguing by induction we claim such an equivalence but, with the notation there, if  the group $\scriptstyle U$  is {\cdt normal} in {\cdy  F} then the induction argument cannot be applied!}}
$$\P^{^{\frak X,W}}\cong \T_{L^{^W}_\F}^{^{\frak X^{^W}}}
\eqno  £2.6.1;$$
note that, according to~£2.3.1, we may assume that $W$ does not belong to $\frak X\,.$

\bigskip
\noindent
{\bf Proposition~£2.7.} {\it If $\P^{^\frak X}\!$ is an extendable perfect $\F^{^\frak X}\!\-$locality
then, for any $\F\-$selfcentralizing subgroup $V$ of $P$ {\it fully normalized\/} in~$\F$ such that
$P^{^V}\in \frak X\,,$
$\P^{^{\frak X,V}} \!\!$ is an extendable perfect $\F^{^{\frak X,V}}\!\!\-$locality.\/}

\medskip
\noindent
{\bf Proof:} From our definition we have an {\it $\F^{^V}\-$locality isomorphism\/}
$$\P^{^{\frak X,V}}\cong \T_{L^{^V}_\F}^{^{\frak X^{^V}}}
\eqno  £2.7.1,$$
which  determines an {\it $N_{\F^{^V}} (W)\-$locality isomorphism\/}
$$N_{\P^{^{\frak X,V}}} (W)\cong N_{\T_{\!L^{^V}_\F}^{^{\frak X^{^V}}}} (W) = \T_{N_{\!L^{^V}_\F} (W)}^{^{\frak X^{^{V,W}}}}
\eqno  £2.7.2\phantom{.}$$
where we identify $P^{^V}$ with its image in $L^{^V}_\F$ and,  for any   $\F^{^V}\-$selfcentralizing subgroup $W$ of~$P^{^V}\!$ {\it fully normalized\/} in~$\F^{^V}$ such that 
$N_{P^{^V}} (W)\in \frak X^{^V}\,,$ we denote by $\frak X^{^{V,W}}$ the set of subgroups of $N_{P^{^V}} (W)$ belonging to $\frak X^{^V}$

\smallskip
But,  it is not difficult to check that the normalizer $N_{\!L^{^V}_\F} (W)\,,$ endowed with the  group homomorphisms
$$N_{P^{^V}} (W)\too N_{L^{^V}_\F} (W)\qq N_{L^{^V}_\F} (W)\too \F^{^V} \!(W)
\eqno £2.7.3\phantom{.}$$
\eject
\noindent
induced by the structural group homomorphisms of $L^{^V}_\F\,,$ is the {\it $\F^{^V}\!\-$localizer\/}
of $W\,.$ We are done.

\bigskip
\noindent
{\bf £3. A reduction procedure\/}

\medskip
£3.1. With the notation above, recall that a {\it basic $P\times P\-$set\/} [3,~21,4] is a finite nonempty $P\times P\-$set $\Omega$ such that $\{1\}\times P$ acts {\it freely\/} on $\Omega\,,$ that we have
$$\Omega^\circ \cong \Omega\qq \vert\Omega\vert/\vert P\vert \not\equiv 0 \bmod{p}
\eqno £3.1.1\phantom{.}$$
where we  denote by $\Omega^\circ$ the $P\times P\-$set obtained by exchanging both factors,
and that, for any subgroup $Q$ of $P$ and any injective  group homomorphism $\varphi\,\colon Q\to P$ such that $\Omega$ contains a $P\times P\-$subset isomorphic to $(P\times P)/\Delta_\varphi (Q)$ where we set 
$\Delta_\varphi (Q) = \{(\varphi (u),u)\}_{u\in Q}\,,$ we have a $Q\times P\-$set isomorphism
$${\rm Res}_{\varphi\times {\rm id}_P} (\Omega)\cong {\rm Res}_{\iota_Q^P\times {\rm id}_P} (\Omega)
\eqno £3.1.2\phantom{.}$$

\medskip 
£3.2. Denoting by $G^{^\Omega}$ the group of  automorphisms of the $\{1\}\times P\-$set
${\rm Res}_{\{1\}\times P}(\Omega)\,,$ it is clear that we have an injective~map from $P\times \{1\}$ 
in~$G^{^\Omega}\,;$  we identify its image with the $p\-$group $P$ so that, from now on, $P$ is contained in $G^{^\Omega}$ and acts freely on $\Omega\,.$ Recall that the {\it full\/} subcategory of the
{\it $G^{^\Omega}\-$transporter\/} over the set of subgroups of $P$ induces a 
{\it Frobenius $P\-$category\/} [3,~Proposition~21.9] and we say that $\Omega$ is an 
{\it $\F\-$basic $P\times P\-$set\/} if, for any pair of subgroups $Q$ and $R$ of $P\,,$ 
we have 
$$\T_{G^{^\Omega}} (Q,R)/C_{G^{^\Omega}}(R)\cong \F (Q,R)
\eqno £3.2.1.$$

\medskip
£3.3. Actually, it follows from [3,~Proposition~21.12] that an {\it $\F\-$basic $P\times P\-$set\/}
always exists; more precisely, we say that an {\it $\F\-$basic $P\times P\-$set\/} $\Omega$ is {\it natural\/} if it fulfills [5,~3.5]
$$\vert \Omega^{\Delta_\varphi (Q)}\vert = \vert Z (Q)\vert
\eqno £3.3.1\phantom{.}$$
 for any $\F\-$selfcentralizing subgroup $Q$ of $P$ and any $\varphi\in \F (P,Q)\,,$ and if it is {\it thick\/}
 [3,~21.7] outside of the set  of $\F\-$selfcentralizing subgroups of $P$ --- namely  the multiplicity of  $(P\times P)/\Delta_\psi (R)$  is at least {\it two\/} if $R$ is not $\F\-$selfcentralzing and $\psi$ belongs to $\F (P,R)\,.$ The existence of {\it natural $\F\-$basic $P\times P\-$sets\/} follows from [5,~Proposition~3.4] together with [3,~Proposition~21.12];
here, we are interested in the following form of [5,~Proposition~3.7]

\bigskip
\noindent
{\bf Proposition~£3.4} {\it Let $\Omega$ be a natural $\F\-$basic $P\times P\-$set, $Q$ and $T$ a pair of   $\F\-$selfcentralizing subgroups of $P$  and $\eta$ an element of~$\F (Q,T).$ The multiplicity of 
$(Q\times P)/\Delta_\eta (T)$ in ${\rm Res}_{Q\times P}(\Omega)$ is at most one, 
and if it is one then we have 
$${\rm Aut}_{Q\times P}\big((Q\times P)/\Delta_\eta (T)\big)\cong Z(T)
\eqno £3.4.1.$$\/}

\medskip
£3.5. From now on, $\Omega$ is  a {\it natural $\F\-$basic $P\times P\-$set\/}. For any subgroup $Q$ \hskip 1pt of \hskip 1pt $P\,,$ \hskip 1pt  it \hskip 1pt is \hskip 1pt clear \hskip 1pt that \hskip 1pt $C_{G^{^\Omega}}(Q)$ \hskip 1pt is \hskip 1pt just \hskip 1pt the \hskip 1pt group \hskip 1pt of \hskip 1pt  automorphisms \hskip 1pt of \hskip 1pt the
\eject
\noindent
$Q\times P\-$set ${\rm Res}_{Q\times P} (\Omega)$ and it is clear that the correspondence 
sending~$Q$ to~$C_{G^{^\Omega}}(Q)$ induces a {\it contravariant\/} functor $C_{G^{^\Omega}}$ from $\F$ to the category $\Gr$ of finite groups. Let us denote by $C^{^{\rm nsc}}_{G^{^\Omega}}(Q)$ the subgroup of 
elements $f\in C_{G^{^\Omega}}(Q)$ which act trivially
on all the $Q\times P\-$orbits of~$\Omega$ isomorphic to $(Q\times P)/\Delta_\eta (T)$ where $T$ is
 $\F\-$selfcentralizing; in particular, if $Q$ is not $\F\-$selfcentralizing then we have
$C^{^{\rm nsc}}_{G^{^\Omega}}(Q) = C_{G^{^\Omega}}(Q)\,;$ in any case, 
$C^{^{\rm nsc}}_{G^{^\Omega}}(Q)$ is normal in $C_{G^{^\Omega}}(Q)$ and, according to
Proposition~£3.4, the quotient $C_{G^{^\Omega}}(Q)/C^{^{\rm nsc}}_{G^{^\Omega}}(Q)$ is Abelian.

\medskip
£3.6. More generally, for any $Q\in \frak X$  denote by $C^{^{\C\frak X}}_{G^{^\Omega}}(Q)$ the subgroup of  
elements $f\in C_{G^{^\Omega}}(Q)$ which act trivially on all the $Q\times P\-$orbits 
of~$\Omega$ isomorphic to $(Q\times P)/\Delta_\eta (T)$ where $T$ belongs to $\frak X\,;$
it is easily checked that the correspondence sending $Q\in \frak X$ to $C^{^{\C\frak X}}_{G^{^\Omega}}(Q)$ defines a {\it subfunctor\/} $C^{^{\C\frak X}}_{\!G^{^\Omega}}\,\colon \F^{^\frak X}\to \Gr$ of the restriction
of $C_{G^{^\Omega}}$  to $\F^{^\frak X}\,,$ and we consider the 
{\it quotient $\F^{^\frak X}\!\-$locality\/} $\overline{\T}_{\!G^{^\Omega}}^{^\frak X} = 
\T_{G^{^\Omega}}^{^\frak X}/C^{^{\C\frak X}}_{\!G^{^\Omega}}$  --- noted 
$\bar\L^{^{\rm n,\frak X}}$ in [5,~5.1.2] --- sending any pair of groups $Q$ and $R$ in $\frak X$ to
$$\overline{\T}_{\!G^{^\Omega}}^{^\frak X} (Q,R) = 
\T_{G^{^\Omega}}^{^\frak X} (Q,R)/C^{^{\C\frak X}}_{\!G^{^\Omega}} (R)
\eqno £3.6.1;$$
 here we are interested in the following form of [5,~Corollaries~5.20 and~5.21].

\bigskip
\noindent
{\bf Proposition~£3.7.} {\it For any perfect  $\F^{^{\frak X}}\!\-$locality $\P^{^\frak X}$ there is a unique naturally $\F^{^{\frak X}}\!\-$isomorphic  class of faithful $\F^{^\frak X}\-$locality functors 
$\lambda\!^{^\frak X}\,\colon \P^{^\frak X}\to \overline{\T}^{^{\frak X}}_{G^{^\Omega}}\,.$ 
Moreover, if $\P'^{^{\frak X}}\!$ is a perfect  $\F^{^{\frak X}}\!\-$locality which is 
$\F^{^{\frak X}}\!\-$locality isomorphic to $\P^{^\frak X}$ then there is a commutative diagram of
 $\F^{^{\frak X}}\!\-$locality functors
$$\matrix{\P^{^{\frak X}}&\buildrel \rho^{_{\frak X}}\over
\cong &\P'^{^{\frak X}}\cr
{\atop \lambda^{^{\!\frak X}}}\searrow&
&\hskip-10pt\swarrow{\atop \lambda'^{^\frak X}}\cr
&\overline{\T}^{^{\frak X}}_{\!G^{^\Omega}}\cr}
\eqno £3.7.1.$$\/}

\medskip
£3.8. With the notation in~£2.2 above, for any $\F\-$selfcentralizing subgroup $W$ of $P$ {\it fully normalized\/} in $\F$ such that $P^{^W}\!\in \frak X\,,$ it follows from [3,~Proposition~21.11] that the subset of $\Omega$
$$\Omega_{_W} = \bigcup_{\chi\in \F(W)} \Omega^{\Delta_\chi (W)}
\eqno £3.8.1\phantom{.}$$
is actually an {\it $\F^{^W}\!\-$basic $P^{^W}\!\!\times P^{^W}\!\!\-$set\/}; {\it mutatis mutandi\/},  denote by $G^{^{\Omega_{_W}}}$ the group of $\{1\}\times P^{^W}\!\!\-$set automorphisms of 
${\rm Res}_{\{1\}\times P^{^W}\!}(\Omega_{_W})$ and identify $P^{^W}\!$ with 
$P^{^W}\!\!\times \{1\}\,;$   since the quotient $N_{G^{^\Omega}} (W)/C_{G^{^\Omega}} (W)$ is isomorphic to $\F (W)$
 (cf.~£3.2.1), it is clear that $N_{G^{^\Omega}} (W)$ stabilizes $\Omega_{_W}$ and therefore we have a canonical group homomorphism from $N_{G^{^\Omega}} (W)$ to $G^{^{\Omega_{_W}}}\,;$ again,  we are interested in the following form of [5,~Proposition~6.15].
 \eject

\bigskip
\noindent
{\bf Proposition~£3.9.} {\it With the notation above, for any pair of  subgroups $Q$ and~$R$ of 
$P^{^W}\!$  containing~$W$  and any element $\varphi$ in $\F^{^W}\! (Q,R)\,,$ there exists at most one $Q\times P^{^W}\!\!\-$orbit in $\Omega_{_W}$ isomorphic to 
$(Q\times P^{^W})/\Delta_\varphi (R)\,,$  $\Omega_{_W}$ is a natural 
$\F^{^W}\!\-$basic $P^{^W}\!\times P^{^W}\!\-$set and, in particular, $C_{G^{^{\Omega_{_W}}}} ( Q)$ is an Abelian $p\-$group.\/}

\medskip
£3.10. It follows from this proposition that, as in~£3.6 above, if $P^{^W}$ belongs to $\frak X$
then we get the {\it quotient $\F^{^{\frak X,W}}\!\-$locality\/} 
$\overline{\T}_{\!G^{^{\Omega_{_W}}}}^{^{\frak X^{^W}}}\,;$ actually, it follows from
 Propositions~£3.4 and~£3.9 above that, with the notation in~£2.2 and~£2.6 above, the canonical group homomorphism from  $N_{G^{^\Omega}} (W)$ to $G^{^{\Omega_{_W}}}$ induces an 
$\F^{^{\frak X,W}}\!\-$locality functor
$$\overline{\frak g}^{^{\frak X,W}}_{_\Omega} : N_{\overline{\T}^{^\frak X}_{\!G^{^\Omega}}\!} (W) \too 
\overline{\T}^{^{\frak X^{^W}}}_{\! G^{^{\Omega_{_W}}}}
\eqno £3.10.1;$$
note that, according to Proposition~£3.7 above, we have {\it faithful $\F^{^{\frak X,W}}\!\-$locality functors\/} from $\P^{^{\frak X,W}} = N_{\P^{^\frak X}} (W)$ to both $\F^{^{\frak X,W}}\!\-$localities $N_{\overline{\T}^{^\frak X}_{\!G^{^\Omega}}\!} (W)$ and  $\overline{\T}^{^{\frak X^{^W}}}_{\! G^{^{\Omega_{_W}}}}$
and we may assume that they agree with $\overline{\frak g}^{^{\frak X,W}}_{_\Omega}\,.$

\medskip
£3.11.  On the other hand, let $L^{^W}_\F$  be the {\it $\F\-$localizer\/} of $W$ 
[3,~Theorem~18.6]; that is to say, $L^{^W}_\F$ is a finite group endowed with an injective and a surjective group homomorphisms
$$\tau^{_W}_{_{\F}} : P^{^W}\too L^{^W}_\F \qq \pi^{_W}_{_{\F}} : L^{^W}_\F\too \F (W)
\eqno £3.11.1,$$
$\tau^{_W}_{_{\F}} (P^{^W})$ is a Sylow $p\-$subgroup of $L^{^W}_\F\,,$ the composition 
$\pi^{_W}_{_{\F}}\!\circ \tau^{_W}_{_{\F}}\!$ is defined by the conjugation in $\F (W)$ and we also have the exact sequence
$$1\too Z(W)\buildrel \tau^{_W}_{_{\F}} \over{\hbox to 20pt{\rightarrowfill}} L^{^W}_\F\buildrel \pi^{_W}_{_{\F}}\over{\hbox to 20pt{\rightarrowfill}} \F(W)\too 1
\eqno £3.11.2.$$
Below, we restate [5,~Proposition~6.19].

\bigskip
\noindent
{\bf Proposition~£3.12.} {\it With the notation above, there is a unique 
$C_{\! G^{^{\Omega_{_W}}}} (W)\-$ conjugacy class of group homomorphisms
$$\lambda^{^W}_{_{\F}} : L^{^W}_\F\too N_{\! G^{^{\Omega_{_W}}}} (W)
\eqno £3.12.1\phantom{.}$$
compatible with the structural group homomorphisms from $P^{^W}$ and to $\F (W)\,.$\/}

\medskip
£3.13.  As in~£2.6  above,  denote by~$\T_{\!L^{^W}_\F}$ the {\it $\F^{^W}\!\-$locality\/} determined by 
$\tau^{_W}_{_{\F}}$ and by the {\it transporter\/} of the group  $L^{^W}_\F\,;$ it is clear that any group homomorphism $\lambda^{^W}_{_{\F}}\,\colon  L^{^W}_\F\to N_{\! G^{^{\Omega_{_W}}}} (W)$ in~£3.12.1 above determines an {\it $\F^{^W}\!\-$locality functor\/}
$$\frak l^{^W}_{_{\F}} : \T_{\!L^{^W}_\F}\too \T_{\! G^{^{\Omega_{_W}}}}
\eqno £3.13.1\phantom{.}$$
\eject
\noindent
and two of them are {\it naturally $\F^{^W}\!\-$isomorphic\/} [5,~2.9]; moreover, if 
$P^{^W}\in \frak X\,,$ it is not difficult to see that the {\it full\/} subcategory 
$\T_{\! L^{^W}_\F}^{^{\frak X^{^W}}}$ of $\T_{\! L^{^W}_\F}$ over $\frak X^{^W}$ as the set of objects is a 
{\it perfect $\F^{^{\frak X,W}}\!\-$locality\/}, and from~£3.13.1 we get an 
{\it $\F^{^{\frak X,W}}\!\-$locality functor\/}
$$\frak l_{_{\F}}^{^{\frak X,W}} : \T_{\! L^{^W}_\F}^{^{\frak X^{^W}}}\too 
\overline{\T}_{\! G^{^{\Omega_{_W}}}}^{^{\frak X^{^W}}}
\eqno £3.13.2.$$

\bigskip
\noindent
{\bf £4. Existence and uniqueness of an extendable perfect $\F^{^{\frak X}}\!\-$locality\/}

\medskip
£4.1. With the notation in~£1.3 above, our main purpose  is to prove that 
\bigskip
\noindent
{\bf Theorem.} {\it There exists an
extendable perfect $\F^{^{\frak X}}\!\-$locality $\P^{^\frak X}\!$, which is unique up to 
$\F^{^{\frak X}}\!\-$locality isomorphisms\/}.
\medskip
\noindent
The existence and the uniqueness of the {\it $\F\-$localizer\/} $L^{^P}_\F$ of $P$ [3,~Theorem~18.6]
proves the existence and the uniqueness of the  {\it extendable perfect $\F^{^\frak X}\-$locality\/} whenever $\frak X = \{P\}\,;$ indeed, $L^{^P}_\F$ is actually a semidirect product $P\!\rtimes\! K$ where
$K\cong \F (P)/\F_P (P)$ is a $p'\-$group and, for any $\F\-$selfcentralizing normal subgroup $W$ of 
$P\,,$ the $\F^{^{\frak X,W}}\!\-$locality equivalence~£2.6.1 is obvious.

\medskip
£4.2. Thus, we may assume that $\frak X\not= \{P\}$ and will argue by induction 
on~$\vert\frak X\vert\,.$ Choose a minimal element $U$ in $\frak X$ {\it fully normalized\/} 
in $\F$ and set 
$$\frak Y = \frak X - \{\theta(U)\mid \theta\in \F(P,U)\}
\eqno £4.2.1;$$
then, by the induction hypothesis, we may assume that  there exists an
{\it extendable perfect $\F^{^{\frak Y}}\-$locality\/}~$\P^{^{\frak Y}}\,,$ endowed with the
structural functors
$$\tau^{_\frak Y} : \T_P^{^\frak Y}\too \P^{^{\frak Y}}\qq
\pi^{_\frak Y} : \P^{^{\frak Y}}\too \F^{^{\frak Y}}
\eqno £4.2.2,$$
which is unique up to $\F^{^\frak Y}\-$locality isomorphisms. At this point, according to 
Proposition~£3.7 above, we may assume that  $\P^{^{\frak Y}}$  is an 
{\it $\F^{^{\frak Y}}\!\-$sublocality\/} of the {\it $\F^{^{\frak Y}}\!\-$locality\/} 
$\overline{\T}^{^\frak Y}_{\!G^{^\Omega}}$ introduced in~£3.6 above; then, denoting by 
$(\overline{\T}^{^\frak X}_{\!G^{^\Omega}}\!)^{^\frak Y}$ the {\it full\/} subcategory of 
$\overline{\T}^{^\frak X}_{\!G^{^\Omega}}\!$ over $\frak Y$ as the set of objects, we  have an obvious functor 
$(\overline{\T}^{^\frak X}_{\!G^{^\Omega}}\!)^{^\frak Y}\too \overline{\T}^{^\frak Y}_{\!G^{^\Omega}}$ and we look to the {\it pull-back\/}
$$\matrix{\P^{^{\frak Y}}&\i &\overline{\T}^{^\frak Y}_{\!G^{^\Omega}}\!\cr
\hskip-5pt\uparrow&\phantom{\big\uparrow}&\hskip-5pt\uparrow\cr
\M^{^{\Omega,\frak Y}}\!&\i &(\overline{\T}^{^\frak X}_{\!G^{^\Omega}}\!)^{^\frak Y}\cr}
\eqno £4.2.3,$$
which defines a  {\it coherent $\F^{^{\frak Y}}\!\-$locality\/} $\M^{^{\Omega,\frak Y}}\!$
[3,~17.9] endowed with obvious structural functors
$$\upsilon^{_{\Omega,\frak Y}} : \T_P^{^\frak Y}\too \M^{^{\Omega,\frak Y}}\qq
\rho^{_{\Omega,\frak Y}} : \M^{^{\Omega,\frak Y}}\too \F^{^{\frak Y}}
\eqno £4.2.4.$$
\eject

\medskip
£4.3. We extend $\M^{^{\Omega,\frak Y}}\!$ to  a  {\it coherent 
$\F^{^{\frak X}}\!\-$sublocality\/} $\M^{^{\Omega,\frak X}}\!$ of 
$\overline{\T}^{^\frak X}_{\!G^{^\Omega}}\!$ which contains $\M^{^{\Omega,\frak Y}}\!$ 
as a {\it full\/} subcategory over $\frak Y$ and fulfills
$$\M^{^{\Omega,\frak X}}\! (Q,V) = \overline{\T}^{^\frak X}_{\!G^{^\Omega}} (Q,V)
\eqno £4.3.1\phantom{.}$$
for any $Q\in \frak X$ and any $V\in \frak X -\frak Y\,,$ and denote by 
$$\upsilon^{_{\Omega,\frak X}} : \T_P^{^\frak X}\too \M^{^{\Omega,\frak X}}\qq
\rho^{_{\Omega,\frak X}} : \M^{^{\Omega,\frak X}}\too \F^{^{\frak X}}
\eqno £4.3.2\phantom{.}$$
the corresponding structural functors; finally, we consider the {\it quotient 
$\F^{^{\frak X}} \-$lo-cality\/}   $\bar\M^{^{\Omega,\frak X}}\!$ of $\M^{^{\Omega,\frak X}}\!$ defined by
$$\bar\M^{^{\Omega,\frak X}}\!(Q,R) = \M^{^{\Omega,\frak X}}\!(Q,R)\big/\upsilon^{_{\Omega,\frak X}}_{_R} \big(Z (R)\big)
\eqno £4.3.3\phantom{.}$$
for any $Q,R\in \frak X\,,$
together  with the induced natural maps --- denoted by $\bar\upsilon^{^{\Omega,\frak X}}$ and 
$\bar\rho^{^{\Omega,\frak X}}\,.$ Then, the proof of the Theorem above can be reduced to the proof of the following fact, that we prove in the next section
\smallskip
\noindent
£4.3.4.\quad {\it  The structural functor $\bar\rho^{^{\Omega,\frak X}}\!$ admits an  
$\F^{^{\frak X}}\!\-$locality functorial section.\/}

\medskip
£4.4. Let us first prove this reduction. Choose an $\F^{^{\frak X}}\!\-$locality functorial section
$\bar\sigma^{_{\Omega,\frak X}}\,\colon \F^{^{\frak X}}\to  \bar\M^{^{\Omega,\frak X}}\,;$  for any 
pair of groups $Q$ and $R$ in $\frak Y\,,$ we know that (cf.~£2.1)
$$\F^{^{\frak X}}\! (Q,R)\cong \P^{^{\frak Y}}\! (Q,R)\big/\tau^{_\frak Y}_{_R} \big(Z (R)\big)
\eqno £4.4.1\phantom{.}$$
and therefore, denoting by $\P^{^{\Omega,\frak Y}}\! (Q,R)$ the converse image of 
$\bar\sigma^{_{\Omega,\frak X}}_{_{Q,R}}\big(\F^{^{\frak X}}\! (Q,R)\big)$ 
in~$\M^{^{\Omega,\frak X}}(Q,R)\,,$ it is clear that the canonical map
$\M^{^{\Omega,\frak X}}\!(Q,R)\to \P^{^\frak Y}\!(Q,R)$ induces a bijection
$\P^{^{\Omega,\frak Y}}\! (Q,R)\cong \P^{^\frak Y}\! (Q,R)\,;$ that is to say,  looking to the {\it pull-back\/}
$$\matrix{\F^{^{\frak X}}
&\buildrel \bar\sigma^{_{\Omega,\frak X}}\over{\hbox to 25pt{\rightarrowfill}}
&\bar\M^{^{\Omega,\frak X}}\cr
\hskip-5pt\uparrow&\phantom{\big\uparrow}&\hskip-10pt\uparrow\cr
\P^{^{\Omega,\frak X}}\!&{\hbox to 25pt{\rightarrowfill}} &\M^{^{\Omega,\frak X}}\cr}
\eqno £4.4.2\phantom{.}$$
--- which defines a  {\it coherent $\F^{^{\frak X}}\!\-$locality\/} $\P^{^{\Omega,\frak X}}\!$
[3,~17.9] endowed with obvious structural functors
$$\tau^{_{\Omega,\frak X}} : \T_P^{^\frak X}\too \P^{^{\Omega,\frak X}}\qq
\pi^{_{\Omega,\frak X}} : \P^{^{\Omega,\frak X}}\too \F^{^{\frak X}}
\eqno £4.4.3\phantom{.}$$
--- and denoting by $(\P^{^{\Omega,\frak X}}\!)^{^\frak Y}$ the {\it full\/} subcategory of
$\P^{^{\Omega,\frak X}}\!$ over $\frak Y$ as the set of objects, it follows from those 
bijections above that we have an {\it $\F^{^\frak Y}\!\-$locality isomorphism\/} 
$(\P^{^{\Omega,\frak X}}\!)^{^\frak Y}\cong \P^{^\frak Y}\,.$ 

\medskip
£4.5. That is to say,  for any $Q\in \frak Y$ {\it fully normalized\/} in $\F\,,$ we already know that 
$\P^{^{\Omega,\frak X}}\!(Q) $ is an {\it $\F\-$localizer\/} of~$Q$ and, 
for  any $V\in \frak X -\frak Y\,,$ it follows from the {\it pull-back\/}~£4.4.2 above that we have the exact sequence
$$1\too Z(V)\too \P^{^{\Omega,\frak X}}\!(V)\too \F(V)\too 1
\eqno £4.5.1\phantom{.}$$\break
\eject
\noindent
and it is easily checked that the group $\P^{^{\Omega,\frak X}}\!(V)\,,$ endowed with the group homomorphisms
$$\tau^{_{\Omega,\frak X}}_{_{V}} : N_P (V)\too \P^{^{\Omega,\frak X}}\!(V) \qq
\pi^{_{\Omega,\frak X}}_{_{V}} : \P^{^{\Omega,\frak X}}\!(V)\too \F (V)
\eqno £4.5.2\phantom{.}$$
determined by the functors $\tau^{_{\Omega,\frak X}}$ and $\pi^{_{\Omega,\frak X}}\,,$
is actually an {\it $\F^{^{\frak X}}\-$localizer\/} of $V$ whenever $V$ is {\it fully normalized\/} in $\F\,;$ 
consequently, it follows from £2.1 above that $\P^{^{\Omega,\frak X}}\!$ is a {\it perfect $\F^{^{\frak X}}\!\-$locality\/}.

\medskip
£4.6. We  claim that $\P^{^{\Omega,\frak X}}\!$ is actually an {\it extendable perfect  
$\F^{^{\frak X}}\!\-$locality\/}; indeed, let $W$ be an $\F\-$selfcentralizing subgroup of $P$ {\it fully 
normalized\/} in $\F$ such that $P^{^W} = N_P(W)$ belongs to $\frak X\,;$ thus, if $P^{^W}$ does not belong to $\frak Y$ then  we have $\frak X^{^W} = \{P^{^W}\}$ and $P^{^W}$ is the unique object in both $\F^{^{\frak X,W}}\!\-$localities $N_{\P^{^{\Omega,\frak X}}} (W)$ and 
$\T_{L^{^W}_\F}^{^{\frak X^{^W}}}\,;$ in this case, since
$$\big(N_{\P^{^{\Omega,\frak X}}} (W)\big) (P^{^W}) \cong P^{^W}\!\!\rtimes K \cong \T_{L^{^W}_\F}^{^{\frak X^{^W}}} (P^{^W})
\eqno £4.6.1\phantom{.}$$
where $K\cong \F^{^{\frak X,W}}\!\! (P^{^W})/\F_{\!P^{^W}}(P^{^W})\,,$ it is  clear that we get the equivalence~£2.6.1.
 Otherwise $\frak Y^{^W}$ is not empty and, setting $\P^{^{\Omega,\frak X,W}} = N_{\P^{^{\Omega,\frak X}}} (W)$ and denoting by 
$\P^{^{\Omega,\frak Y,W}}\!$ the {\it full\/} subcategory of~$\P^{^{\Omega,\frak X,W}} \!$ over 
$\frak Y^{^W}\,,$ from~£4.4 above we get an {\it $\F^{^{\frak Y,W}}\!\-$locality isomorphism\/} 
$$\P^{^{\Omega,\frak Y,W}}\!\cong N_{\P^{^{\frak Y}}\!} (W)
\eqno £4.6.2;$$
but, since $\P^{^\frak Y}$ is {\it extendable\/}, it follows from our definition in~£2.6 above that we still get 
an {\it $\F^{^{\frak Y,W}}\!\-$locality isomorphism\/} 
$$N_{\P^{^{\frak Y}}\!} (W)\cong  \T_{\!L^{^W}_\F}^{^{\frak Y^{^W}}}
\eqno £4.6.3.$$

\par
£4.7. Always assuming that $\frak Y^{^W}$ is not empty, note that in~£3.10 
above~$\overline{\frak g}^{^{\frak Y,W}}_{_\Omega}$ sends $N_{\P^{^{\frak Y}}\!} (W)\!$ isomorphically to its image in $\overline{\T}^{^{\frak Y^{^W}}}_{G^{^{\Omega_{_W}}}} $ --- still noted 
$\P^{^{\Omega,\frak Y,W}}\,;$
then, from this inclusion, {\it mutatis mutandi\/}  we can define a  {\it coherent~$\F^{^{\frak Y,W}}\!\-$ locality\/} $\M^{^{\Omega,\frak Y,W}}\!$ as in~£4.2.3, and 
{\it coherent~$\F^{^{\frak X,W}}\!\-$localities\/} $\M^{^{\Omega,\frak X,W}}\!\i 
\overline{\T}^{^{\frak X^{^W}}}_{\!G^{^{\Omega_{_W}}}}$ and~$\bar\M^{^{\Omega,\frak X,W}}\!$ as in~£4.3; moreover, it is clear that $\bar\sigma^{_{\Omega,\frak X}}$ induces an 
$\F^{^{\frak X,W}}\!\-$locality functorial section 
$\bar\sigma^{_{\Omega,\frak X,W}}\,\colon \F^{^{\frak X,W}}\to  \bar\M^{^{\Omega,\frak X,W}}$ and that we can define a {\it coherent $\F^{^{\frak X,W}}\!\-$locality\/} $\P^{^{\Omega,\frak X,W}}$ as in~£4.4.2 above which  still fulfills 
$$(\P^{^{\Omega,\frak X,W}})^{^{\frak Y^{^W}}}\cong \P^{^{\Omega,\frak Y,W}}
\eqno £4.7.1;$$
we  denote by $\tau^{_{\Omega,\frak X,W}}\,\colon\T_{P^{^W}}^{^{\frak X^{^W}}}\to 
\P^{^{\Omega,\frak X,W}}$ and by
$\pi^{_{\Omega,\frak X,W}}\,\colon \P^{^{\Omega,\frak X,W}}\to \F^{^{\frak X,W}}$
the structural functors. 
\eject

\medskip
£4.8. On the other hand, since $\T_{\!L^{^W}_\F}^{^{\frak X^{^W}}}$ is a {\it perfect  
$\F^{^{\frak X,W}}\!\-$locality\/} (cf.~£3.13), it follows from Proposition~£3.7 (or from~£3.13.2) that 
$\T_{\! L^{^W}_\F}^{^{\frak X^{^W}}}\!$ is actually an {\it $\F^{^{\frak X,W}}\!\!\-$sublocality\/} of 
$\overline{\T}^{^{\frak X^{^W}}}_{\!G^{^{\Omega_{_W}}}}\,;$ in particular, denoting by 
$(\overline{\T}^{^{\frak X^{^W}}}_{\! G^{^{\Omega_{_W}}}})^{^{\frak Y^{^W}}}$ and by
$(\T_{\! L^{^W}_\F}^{^{\frak X^{^W}}})^{^{\frak Y^{^W}}}$ the respective {\it full\/} subcategories
of $\overline{\T}^{^{\frak X^{^W}}}_{\!G^{^{\Omega_{_W}}}}$ and of $\T_{\! L^{^W}_\F}^{^{\frak X^{^W}}}$ over $\frak Y^{^W}\,,$ it is easily checked that the canonical functor
$$(\overline{\T}^{^{\frak X^{^W}}}_{\! G^{^{\Omega_{_W}}}})^{^{\frak Y^{^W}}}\!\!\too 
\overline{\T}^{^{\frak Y^{^W}}}_{\! G^{^{\Omega_{_W}}}}
\eqno £4.8.1\phantom{.}$$
sends $(\T_{\! L^{^W}_\F}^{^{\frak X^{^W}}})^{^{\frak Y^{^W}}}$ isomorphically onto 
$\T_{\! L^{^W}_\F}^{^{\frak Y^{^W}}}\i \overline{\T}^{^{\frak Y^{^W}}}_{\! G^{^{\Omega_{_W}}}}\,.$

\medskip
£4.9. Moreover, from £4.4 we know that the canonical functor
$$(\overline{\T}^{^{\frak X}}_{\!G^{^{\Omega}}})^{^{\frak Y}}\!\!\too 
\overline{\T}^{^{\frak Y}}_{\!G^{^{\Omega}}}
\eqno £4.9.1\phantom{.}$$
sends $(\P^{^{\Omega,\frak X}}\!)^{^\frak Y}$ isomorphically onto $\P^{^\frak Y}\,;$ 
but, it follows from our definition in~£3.10  that, denoting by 
$(\overline{\frak g}^{^{\frak X,W}}_{_\Omega})^{^\frak Y}$ the restriction of 
$\overline{\frak g}^{^{\frak X,W}}_{_\Omega}$ to the normalizer in $(\overline{\T}^{^{\frak X}}_{G^{^{\Omega}}})^{^{\frak Y}}\!$ of~$W\,,$ we have a commutative diagram of functors
$$\matrix{N_{\overline{\T}^{^\frak Y}_{\!G^{^\Omega}}\!} (W)
&\buildrel \overline{\frak g}^{^{\frak Y,W}}_{_\Omega}\over{\hbox to 30pt{\rightarrowfill}}
&\overline{\T}^{^{\frak Y^{^W}}}_{\! G^{^{\Omega_{_W}}}}\cr
\big\uparrow&\phantom{\Big\uparrow}&\big\uparrow\cr
N_{(\overline{\T}^{^\frak X}_{\!G^{^\Omega}}\!)^{^\frak Y}} (W)
&\buildrel (\overline{\frak g}^{^{\frak X,W}}_{_\Omega})^{^\frak Y}\over{\hbox to 30pt{\rightarrowfill}} 
&(\overline{\T}^{^{\frak X^{^W}}}_{\! G^{^{\Omega_{_W}}}})^{^{\frak Y^{^W}}}\cr}
\eqno £4.9.2,$$
where the vertical arrows are defined by the functors~£4.8.1 and~£4.9.1; hence, since the functor~£4.9.1 sends $(\P^{^{\Omega,\frak X}}\!)^{^\frak Y}$ isomorphically onto $\P^{^\frak Y}$
(cf.~£4.4), this functor sends $N_{(\P^{^{\Omega,\frak X}}\!)^{^\frak Y}} (W)$ isomorphically onto 
$N_{\P^{^\frak Y}} (W)$ and we already know that $\overline{\frak g}^{^{\frak Y,W}}_{_\Omega}$  sends $N_{\P^{^\frak Y}} (W)$  isomorphically onto $\P^{^{\Omega,\frak Y,W}}\!\!$ (cf.~£4.7), which is isomorphic to 
$\T_{\! L^{^W}_\F}^{^{\frak Y^{^W}}}$ (cf.~£4.6.2).

\medskip
£4.10. At this point, it follows from Proposition~£3.7 that there exist  an 
{\it $\F^{^{\frak Y,W}}\!\-$locality functor\/} $\frak l_{_{\F}}^{^{\frak Y,W}}\,\colon 
\T_{\! L^{^W}_\F}^{^{\frak Y^{^W}}}\to \overline{\T}^{^{\frak Y^{^W}}}_{\! G^{^{\Omega_{_W}}}}$
which sends $\T_{\! L^{^W}_\F}^{^{\frak Y^{^W}}}$ isomorphically to~$\P^{^{\Omega,\frak Y,W}}\,,$ 
and that  this functor  is {\it naturally $\F^{^{\frak Y,W}}\!\-$isomorphic\/} to the inclusion 
$\T_{\! L^{^W}_\F}^{^{\frak Y^{^W}}}\i \overline{\T}^{^{\frak Y^{^W}}}_{\! G^{^{\Omega_{_W}}}}$ 
in~£4.8 above; that is to say, according to our definition in~[5,~2.9] and since the kernel of the structural group homomorphism from\break
\eject 
\noindent
 $\overline{\T}^{^{\frak Y^{^W}}}_{\! G^{^{\Omega_{_W}}}}  (P^{^W})$ to  $\F^{^{\frak Y,W}}\! (P^{^W})$ is the image of $C_{\! G^{^{\Omega_{_W}}}} (P^{^W})\i \T_{\! G^{^{\Omega_{_W}}}} (P^{^W})\,,$ there is 
$z\in C_{\! G^{^{\Omega_{_W}}}} (P^{^W})$ such that, denoting 
by~$\overline{z}^{^{\frak Y^{^W}}}_{_Q}$ the image of $z$ in 
$\overline{\T}^{^{\frak Y^{^W}}}_{\! G^{^{\Omega_{_W}}}} (Q)$ for any $Q\in \frak Y^{^W}\,,$ we get
$$\P^{^{\Omega,\frak Y,W}} (Q,R) = \overline{z}^{^{\frak Y^{^W}}}_{_Q}\!\!\.
\T_{L^{^W}_\F}^{^{\frak Y^{^W}}} (Q,R)\. (\overline{z}^{^{\frak Y^{^W}}}_{_R}\!)^{-1}
\eqno £4.10.1\phantom{.}$$
in $\overline{\T}^{^{\frak Y^{^W}}}_{\! G^{^{\Omega_{_W}}}} (Q,R)\,,$ for any pair of groups 
$Q$ and $R$ in $\frak Y^{^W}\,.$

\medskip
£4.11. But, we also can consider the images $\overline{z}^{^{\frak X^{^W}}}_{_Q}$ of $z$ in 
$\overline{\T}^{^{\frak X^{^W}}}_{\! G^{^{\Omega_{_W}}}} (Q)$ for any 
$Q\in \frak X^{^W}\,.$ Hence, up to replacing  our choice of 
$\T_{\! L^{^W}_\F}^{^{\frak X^{^W}}}\!$ as a  {\it $\F^{^{\frak X,W}}\!\!\-$sublocality\/} of 
$\overline{\T}^{^{\frak X^{^W}}}_{\!G^{^{\Omega_{_W}}}}$ by the choice 
of~$\overline{z}^{^{\frak X^{^W}}}_{_Q}\!\!\.\T_{\! L^{^W}_\F}^{^{\frak X^{^W}}} (Q,R)\. (\overline{z}^{^{\frak X^{^W}}}_{_R}\!)^{-1}$ in $\overline{\T}^{^{\frak X^{^W}}}_{\! G^{^{\Omega_{_W}}}} (Q,R)\,,$ for any pair of groups $Q$ and $R$ in~$\frak X^{^W}\!\,,$  in 
$\overline{\T}^{^{\frak Y^{^W}}}_{\! G^{^{\Omega_{_W}}}}$ we actually may assume that we get
$$\P^{^{\Omega,\frak Y,W}} = \T_{\! L^{^W}_\F}^{^{\frak Y^{^W}}} 
\eqno £4.11.1.$$
In this situation, it follows from our definitions in~£4.7 above that in 
$\overline{\T}^{^{\frak X^{^W}}}_{\!G^{^{\Omega_{_W}}}}$ the  {\it coherent 
$\F^{^{\frak X,W}}\!\-$sublocality\/} $\M^{^{\Omega,\frak X,W}}$ contains 
$\T_{\! L^{^W}_\F}^{^{\frak X^{^W}}}\!\,.$

\medskip
£4.12. In particular, if $\frak X^{^W}\! = \frak Y^{^W}\!$ then we have
$$\P^{^{\Omega,\frak X,W}} = \P^{^{\Omega,\frak Y,W}} = \T_{\! L^{^W}_\F}^{^{\frak Y^{^W}}} 
= \T_{\! L^{^W}_\F}^{^{\frak X^{^W}}}
\eqno £4.12.1,$$
so that we are done. Assume that $\frak X^{^W}\! \not= \frak Y^{^W}\!\,;$ then, by the very definition  
of~$\overline{\T}^{^{\frak X^{^W}}}_{\!G^{^{\Omega_{_W}}}}$
(cf.~£3.6.1 and~£4.3), for any $V\in \frak X^{^W}\! - \frak Y^{^W}\!$ we have
$${\rm Ker}(\bar\rho^{^{\Omega,\frak X,W}}_{_V})  = \T_{G^{^{\Omega_{_W}}}}^{^{\frak X^{^W}}} (V)/
C^{^{\C\frak X^{^W}}}_{\!G^{^{\Omega_{_W}}}} (V)   =  \prod_{\tilde\theta\in \tilde\F^{^W}(P^{^W},V)}Z(V)
\eqno £4.12.2\phantom{.}$$
and therefore, since $p$ does not divide $\vert \tilde\F^{^W}\!(P^{^W}\!,V)\vert$ 
[3,~Proposition~6.7], we have  a surjective group homomorphism
 $$\nabla^{^{\Omega,\frak X,W}}_V : {\rm Ker}(\bar\rho^{^{\Omega,\frak X,W}}_{_V})\too Z(V)
 \eqno £4.12.3\phantom{.}$$
mapping $z = (z_{\tilde\theta})_{\tilde\theta\in \tilde\F^{^W}\!(P^{^W}\!,V)}$ on 
$$\nabla^{^{\Omega,\frak X,W}}_V (z) = {1\over\vert \tilde\F^{^W}\!(P^{^W}\!,V)\vert}\.
\sum_{\tilde\theta\in \tilde\F^{^W}\!(P^{^W}\!,V)} z_{\tilde\theta}
\eqno £4.12.4.$$
\eject

\medskip
£4.13. At this point, considering  the {\it contravariant  functor\/}  
$$\frak d^{^{\Omega,\frak X,W}} : \overline{\T}^{^{\frak X^{^W}}}_{\!G^{^{\Omega_{_W}}}}\to \Ab
\eqno £4.13.1$$
mapping any $Q\in \frak Y^{^W}$ on $\{0\}$ and any $V\in \frak X^{^W} -\frak Y^{^W}$ on 
${\rm Ker} (\nabla^{^{\Omega,\frak X,W}}_V)\,,$ and the {\it quotient\/} $\F^{^{\frak X,W}}\! \-$locality  
$\overline{\T}^{^{\frak X^{^W}}}_{\!G^{^{\Omega_{_W}}}}\big/\frak d^{^{\Omega,\frak X,W}}$ [5,~2.10],
it is easily checked that the {\it coherent $\F^{^{\frak X,W}}\!\-$localities\/} $\P^{^{\Omega,\frak X,W}}$
(cf.~£4.7) and $\T_{\! L^{^W}_\F}^{^{\frak X^{^W}}}$ have the same image in this quotient; indeed, 
it follows from equalities~£4.11.1 above that their images coincide over $\frak Y^{^W}$ and, since for any
$V\in \frak X^{^W} -\frak Y^{^W}$ we have
$$\big(\overline{\T}^{^{\frak X^{^W}}}_{\!G^{^{\Omega_{_W}}}}\big/\frak d^{^{\Omega,\frak X,W}}\big)(V)
\cong L^{^V}_{\F^{^W}}
\eqno £4.13.2,$$
$\P^{^{\Omega,\frak X,W}} (V)$ and $\T_{\! L^{^W}_\F}^{^{\frak X^{^W}}} (V)$ map both isomorphically onto  
$\big(\overline{\T}^{^{\frak X^{^W}}}_{\!G^{^{\Omega_{_W}}}}\big/\frak d^{^{\Omega,\frak X,W}}\big)(V)\,.$
In particular, we get $\P^{^{\Omega,\frak X,W}}\cong \T_{\! L^{^W}_\F}^{^{\frak X^{^W}}}$ since the functors 
from $\P^{^{\Omega,\frak X,W}}$ and $\T_{\! L^{^W}_\F}^{^{\frak X^{^W}}}$ to the quotient $\overline{\T}^{^{\frak X^{^W}}}_{\!G^{^{\Omega_{_W}}}}\big/\frak d^{^{\Omega,\frak X,W}}$ are faithful. This  proves our claim in~£4.6.

\medskip
£4.14. It remains to prove the uniqueness; thus, assume that $\P^{^\frak X}$ and~$\P'^{^\frak X}$
are two {\it extendable perfect $\F^{^\frak X}\-$localities\/}; it follows from Proposition~£3.7 that 
we may assume that both are $\F^{^\frak X}\-$sublocalities of the $\F^{^\frak X}\-$locality 
$\overline{\T}_{\!G^{^\Omega}}^{^\frak X}$ introduced in~£3.6 above. On the other hand, since the 
respective {\it full\/} sub-categories~$\P^{^\frak Y}$ of $\P^{^\frak X}$ and~$\P'^{^\frak Y}$ of 
$\P'^{^\frak X}$ over $\frak Y$ as the set of objects are still  two {\it extendable perfect 
$\F^{^\frak Y}\-$localities\/}, it follows from our induction hypothesis that they are 
{\it $\F^{^\frak Y}\-$locality isomorphic\/}. Consequently, considering the inclusions of $\P^{^\frak Y}$ and~$\P'^{^\frak Y}$ in $\overline{\T}_{\!G^{^\Omega}}^{^\frak Y}$ induced by the inclusions
$$\P^{^\frak Y} = (\P^{^\frak X})^{^\frak Y}\i (\overline{\T}_{\!G^{^\Omega}}^{^\frak X})^{^\frak Y}
\j (\P'^{^\frak X})^{^\frak Y} = \P'^{^\frak Y}
\eqno £4.14.1$$
and by  the canonical functor $(\overline{\T}^{^{\frak X}}_{\!G^{^{\Omega}}})^{^{\frak Y}}\!\!\to
\overline{\T}^{^{\frak Y}}_{\!G^{^{\Omega}}}$ (cf.~£4.9.1), the existence of an 
{\it $\F^{^\frak Y}\-$locality isomorphism\/} $\P^{^\frak Y}\cong\P'^{^\frak Y}$ determines two 
$\F^{^\frak Y}\-$locality functors from $\P^{^\frak Y}$ to 
$(\overline{\T}_{\!G^{^\Omega}}^{^\frak X})^{^\frak Y}\,;$ then, it follows again from Proposition~£3.7 that the functors ainsi obtained are {\it naturally 
$\F^{^\frak Y}\-$isomorphic\/}.

\medskip
£4.15.  That is to say, as in~£4.10 above,  since the kernel of the structural group homomorphism from 
$\overline{\T}^{^{\frak Y}}_{\! G^{^{\Omega}}}  (P)$ to  $\F^{^{\frak Y}}\! (P)$ is the image of 
$C_{\! G^{^{\Omega}}} (P)$ in~$\T^{^{\frak Y}}_{\! G^{^{\Omega}}} (P)\,,$ there is 
$z\in C_{\! G^{^{\Omega}}} (P)$ such that, denoting by~$\overline{z}^{^{\frak Y}}_{_Q}$ the image of $z$ in 
$\overline{\T}^{^{\frak Y}}_{\! G^{^{\Omega}}} (Q)$ for any $Q\in \frak Y\,,$ in 
$\overline{\T}^{^{\frak Y}}_{\! G^{^{\Omega}}} (Q,R)$  we get
$$\P^{^{\frak Y}} (Q,R) = \overline{z}^{^{\frak Y}}_{_Q}\.
\P'^{^{\frak Y}} (Q,R)\. (\overline{z}^{^{\frak Y}}_{_R}\!)^{-1}
\eqno £4.15.1\phantom{.}$$
for any pair of groups $Q$ and $R$ in $\frak Y\,.$  As above, considering the images 
$\overline{z}^{^{\frak X}}_{_Q}$ of $z$ in $\overline{\T}^{^{\frak X}}_{\! G^{^{\Omega}}} (Q)$ for any  
$Q\in \frak X$ and modifying  our choice of $\P'^{^{\frak X}}\!$ as a  {\it $\F^{^{\frak X}}\!\!\-$sublocality\/} of 
$\overline{\T}^{^{\frak X}}_{\!G^{^{\Omega}}}$ by the choice 
of~$\overline{z}^{^{\frak X}}_{_Q}\.\P'^{^{\frak X}} (Q,R)\. (\overline{z}^{^{\frak X}}_{_R}\!)^{-1}$ in 
$\overline{\T}^{^{\frak X}}_{\! G^{^{\Omega}}} (Q,R)$ for any pair of groups $Q$ and $R$ in~$\frak X\,,$ we actually may assume that in $\overline{\T}^{^{\frak Y^{^W}}}_{\! G^{^{\Omega_{_W}}}}$ we have 
$\P^{^{\frak Y}} = \P'^{^{\frak Y}} \,.$

\medskip
£4.16. Moreover, as in~£4.12 above,  by the very definition of~$\overline{\T}^{^{\frak X}}_{\!G^{^{\Omega}}}$ (cf.~£3.6.1 and~£4.3), for any $V\in \frak X - \frak Y$ we have
$${\rm Ker}(\bar\rho^{^{\Omega,\frak X}}_{_V})  = 
\T_{G^{^{\Omega}}}^{^{\frak X}} (V)/C^{^{\C\frak X}}_{\!G^{^\Omega}} (V)   =  
\prod_{\tilde\theta\in \tilde\F (P,V)}Z(V)
\eqno £4.16.1\phantom{.}$$
and therefore, since $p$ does not divide $\vert \tilde\F (P,V)\vert$ 
[3,~Proposition~6.7], we have  a surjective group homomorphism
 $$\nabla^{^{\Omega,\frak X}}_V : {\rm Ker}(\bar\rho^{^{\Omega,\frak X}}_{_V})\too Z(V)
 \eqno £4.16.2\phantom{.}$$
mapping $z = (z_{\tilde\theta})_{\tilde\theta\in \tilde\F (P,V)}$ on 
$$\nabla^{^{\Omega,\frak X}}_V (z) = {1\over\vert \tilde\F(tP,V)\vert}\.
\sum_{\tilde\theta\in \tilde\F (P,V)} z_{\tilde\theta}
\eqno £4.16.3.$$

\medskip
£4.17. At this point, considering  the {\it contravariant Dirac functor\/}  
$$\frak d^{^{\Omega,\frak X}} : \overline{\T}^{^{\frak X}}_{\!G^{^{\Omega}}}\to \Ab
\eqno £4.17.1$$
mapping any $Q\in \frak Y$ on $\{0\}$ and any $V\in \frak X -\frak Y$ on 
${\rm Ker} (\nabla^{^{\Omega,\frak X}}_V)\,,$ and the~{\it quotient\/} $\F^{^{\frak X}} \-$locality  
$\overline{\T}^{^{\frak X}}_{\!G^{^{\Omega}}}\big/\frak d^{^{\Omega,\frak X}}$ [5,~2.10],
it is easily checked that the {\it coherent $\F^{^{\frak X}}\-$localities\/} $\P^{^{\frak X}}$
 and $\P'^{^{\frak X}}$ have the same image in this quotient; indeed, 
it follows from~£4.15 above that their images coincide over $\frak Y$ and, since for any
$V\in \frak X -\frak Y$ we have
$$\big(\overline{\T}^{^{\frak X}}_{\!G^{^{\Omega}}}\big/\frak d^{^{\Omega,\frak X}}\big)(V)
\cong L^{^V}_\F
\eqno £4.17.2,$$
$\P^{^{\frak X}} (V)$ and $\P'^{^{\frak X}} (V)$ map both isomorphically onto  
$\big(\overline{\T}^{^{\frak X}}_{\!G^{^{\Omega}}}\big/\frak d^{^{\Omega,\frak X}}\big)(V)\,.$
In particular, we get $\P^{^{\frak X}}\cong \P'^{^{\frak X}}$ since the functors 
from $\P^{^{\frak X}}$ and $\P'^{^{\frak X}}$ to the quotient $\overline{\T}^{^{\frak X}}_{\!G^{^{\Omega}}}\big/\frak d^{^{\Omega,\frak X}}$ are faithful. This  proves the uniqueness.
\eject

\bigskip
\noindent
{\bf £5. Existence and uniqueness of the sections from $\F^{^\frak X}$ to 
$\bar\M^{^{\Omega,\frak X}}$\/}

\medskip
£5.1. With the hypothesis and notation in~£4.3 above, our purpose in this section is to prove that
\bigskip
\noindent
{\bf Theorem.} {\it  The structural functor $\bar\rho^{^{\Omega,\frak X}}\!\,\colon 
\bar\M^{^{\Omega,\frak X}}\! \to \F^{^{\frak X}}$ admits an  $\F^{^{\frak X}}\!\-$locality functorial section
$\bar\sigma^{_{\Omega,\frak X}}\,\colon \F^{^{\frak X}}\to  \bar\M^{^{\Omega,\frak X}}\!\,.$\/}

\medskip
\noindent
Actually, since we assume that $U\not= P\,,$ we also have $U\not= N_P(U)  = P^{^U}$ and therefore 
$P^{^U}$ belongs to $\frak Y^{^U}\,;$ thus, this theorem is just the existence part of [5,~Theorem~6.22] but we restate the proof in our new context; indeed, here we assume that $\P^{^\frak Y}$ is an {\it extendable perfect 
$\F^{^\frak Y}\!\-$locality\/} and therefore the {\it $\F^{^{\frak Y,U}}\!\-$locality isomorphism\/} in [5,~6.18]
$$\P^{^{\frak Y,U}}\! = N_{\P^{^\frak Y}} (U) \cong \T_{L^{^U}_\F}^{^{\frak Y^{^U}}}
\eqno £5.1.1\phantom{.}$$
follows from our definition in~£2.6; in particular, as in~£4.11 above, in 
$\overline{\T}^{^{\frak Y^{^U}}}_{\!G^{^{\!\Omega_{_U}}}}\!$ we may assume that 
$\P^{^{\frak Y,U}} = \T_{\! L^{^U}_\F}^{^{\frak Y^{^U}}}\,.$

\medskip
£5.2. Since $\frak Y^{^U}$ is not empty, as in~£4.7 above we can define the 
{\it coherent $\F^{^{\frak Y,U}}\!\-$locality\/} $\M^{^{\Omega,\frak Y,U}}\!$ {\it via\/} the 
{\it pull-back\/} (cf.~£4.2.3)
$$\matrix{\P^{^{\frak Y,U}}&\i &\overline{\T}^{^{\frak Y^{^U}}}_{\!G^{^{\!\Omega_{_U}}}}\!\cr
\hskip-5pt\uparrow&\phantom{\Big\uparrow}&\hskip-5pt\uparrow\cr
\M^{^{\Omega,\frak Y,U}}\!&\i 
&(\overline{\T}^{^{\frak X^{^U}}}_{\!G^{^\Omega}}\!)^{^{\frak Y^{^U}}}\cr}
\eqno £5.2.1\phantom{.}$$
and the {\it coherent $\F^{^{\frak X,U}}\!\-$localities\/} $\M^{^{\Omega,\frak X,U}}\!\i 
\overline{\T}^{^{\frak X^{^U}}}_{\! G^{^{\Omega_{_W}}}}$ and~$\bar\M^{^{\Omega,\frak X,U}}\!$ as in~£4.3, with the second structural functors 
$$\rho^{_{\Omega,\frak X,U}} : \M^{^{\Omega,\frak X,U}}\too \F^{^{\frak X,U}}\qq
\bar\rho^{_{\Omega,\frak X,U}} : \bar\M^{^{\Omega,\frak X,U}}\too \F^{^{\frak X,U}}
\eqno £5.2.2.$$
Now recall that, denoting by $\tilde\F^{^{\frak X}}$ and  $\tilde\F^{^{\frak X,U}}$ the respective
{\it exterior quotients\/} of $\F^{^{\frak X}}$ and  $\F^{^{\frak X,U}}$ [3,~1.3], the {\it coherency\/} of 
$\bar\M^{^{\Omega,\frak X}}\!$ and $\bar\M^{^{\Omega,\frak X,U}}\!$ determines {\it contravariant\/} functors [5,~2.8.3]
$$\Ker (\bar\rho^{^{\Omega,\frak X}}) : \tilde\F^{^{\frak X}}\too \Ab\qq 
\Ker (\bar\rho^{^{\Omega,\frak X,U}}) : \tilde\F^{^{\frak X,U}}\too \Ab
\eqno £5.2.3;$$
as usual, the existence of $\bar\sigma^{_{\Omega,\frak X}}\!$ depends on the vanishing of the {\it cohomology class\/} of a suitable $\Ker (\bar\rho^{^{\Omega,\frak X}})\-$valued $2\-$cocycle and,
from the reduction procedure developed in section~£3, we will move to the corresponding $\Ker (\bar\rho^{^{\Omega,\frak X,U}})\-$va-lued $2\-$cocycle.

\medskip
£5.3. From the commutative diagram~£4.2.3 we get the following commutative diagram of the 
normalizers of $U$
$$\matrix{N_{\P^{^{\frak Y}}}(U)&\i &N_{\overline{\T}^{^\frak Y}_{\!G^{^\Omega}}}(U)\!\cr
\hskip-5pt\uparrow&\phantom{\big\uparrow}&\hskip-5pt\uparrow\cr
N_{\! \M^{^{\Omega,\frak Y}}} (U)\!&\i &N_{(\overline{\T}^{^\frak X}_{\!G^{^\Omega}}\!)^{^\frak Y}} (U)\cr}
\eqno £5.3.1;$$
moreover, we are setting $N_{\P^{^{\frak Y}}}(U) = \P^{^{\frak Y,U}}$ and we have the commutative diagram~£4.9.2 for~$W = U\,.$ Consequently, the $\F^{^{\frak Y,U}}\!\-$ and $\F^{^{\frak X,U}}\!\-$locality functors
 (cf.~£3.10.1)
$$\overline{\frak g}^{^{\frak Y,U}}_{_\Omega} : N_{\overline{\T}^{^\frak Y}_{\!G^{^\Omega}}\!} (U) \too \overline{\T}^{^{\frak Y^{^U}}}_{\! G^{^{\Omega_{_U}}}}\qq
\overline{\frak g}^{^{\frak X,U}}_{_\Omega} : N_{\overline{\T}^{^\frak X}_{\!G^{^\Omega}}\!} (U) \too \overline{\T}^{^{\frak X^{^U}}}_{\! G^{^{\Omega_{_U}}}}
\eqno £5.3.2\phantom{.}$$
successively induce  the new $\F^{^{\frak Y,U}}\-$locality functor (cf.~£5.2.1)
$$\frak h^{^{\frak Y,U}}_{_\Omega} : N_{\!\M^{^{\Omega,\frak Y}}} (U) \too \M^{^{\Omega,\frak Y,U}}
\eqno £5.3.3,$$
and, moreover, the $\F^{^{\frak X,U}}\-$locality functors (cf.~£4.3)
$$\frak h^{^{\frak X,U}}_{_\Omega} : N_{\!\M^{^{\Omega,\frak X}}} (U) \too \M^{^{\Omega,\frak X,U}}\qq \bar\frak h^{^{\frak X,U}}_{_\Omega} : N_{\!\bar\M^{^{\Omega,\frak X}}} (U) \too 
\bar\M^{^{\Omega,\frak X,U}}
\eqno £5.3.4.$$
Similarly, since we are assuming that $\P^{^{\frak Y,U}} = 
\T_{\!L^{^U}_\F}^{^{\frak Y^{^U}}} $ (cf.~£4.11.1), the $\F^{^{\frak X,U}}\!\-$lo-cality functor
(cf.~3.13.2)
$$\frak l_{_{\F}}^{^{\frak X,U}} : \T_{\! L^{^U}_\F}^{^{\frak X^{^U}}}\too 
\overline{\T}_{\! G^{^{\Omega_{_U}}}}^{^{\frak X^{^U}}}
\eqno £5.3.5$$
and the {\it pull-back\/}~£5.2.1 above determine  new $\F^{^{\frak X,U}}\!\-$locality functors (cf.~£4.3)
$$\frak m_{_{\F}}^{^{\frak X,U}} : \T_{\! L^{^U}_\F}^{^{\frak X^{^U}}}\too 
\M^{^{\Omega,\frak X,U}}\qq
\bar\frak m_{_{\F}}^{^{\frak X,U}} : \F^{^{\frak X,U}}\!\too 
\bar\M^{^{\Omega,\frak X,U}}
\eqno £5.3.6.$$

\medskip
£5.4. At this point, denoting by  $\tilde\iota_{_{\frak X,U}}\,\colon \tilde\F^{^{\frak X,U}}\!\to 
\tilde\F^{^{\frak X}}$ the canonical functor, it is well-known that, for any $n\in \Bbb N\,,$ the restriction induces a group homomorphism (cf.~£5.1.2)
$$\Bbb H^n \big(\tilde\F^{^{\frak X}}\!, \Ker (\bar\rho^{^{\Omega,\frak X}})\big)\too 
\Bbb H^n \big(\tilde\F^{^{\frak X,U}}\!, \Ker (\bar\rho^{^{\Omega,\frak X}})
\circ \tilde\iota_{_{\frak X,U}}\big)
\eqno £5.4.1;$$
moreover, $\bar\frak h^{^{\frak X,U}}_{_\Omega}$ induces a {\it natural map\/} [5,~2.10.1]
$$\nu_{\bar\frak h^{^{\frak X,U}}_{_\Omega}} : \Ker (\bar\rho^{^{\Omega,\frak X}})
\circ \tilde\iota_{_{\frak X,U}} \too \Ker (\bar\rho^{^{\Omega,\frak X,U}})
\eqno £5.4.2\phantom{.}$$
and therefore, for any $n\in \Bbb N\,,$ we also get a group homomorphism
$$\Bbb H^n \big(\tilde\F^{^{\frak X,U}}\!, \Ker (\bar\rho^{^{\Omega,\frak X}})
\circ \tilde\iota_{_{\frak X,U}}\big)\too \Bbb H^n \big(\tilde\F^{^{\frak X,U}}\!, 
\Ker (\bar\rho^{^{\Omega,\frak X,U}})\big)
\eqno £5.4.3.$$
\eject
\noindent
In [5,~Proposition~6.9, 6.12.3 and~6.21.7] we prove that, for any $n\in \Bbb N\,,$ the composition 
of the homomorphisms~£5.4.1 and~£5.4.3 determines an isomorphism
$$\Bbb H^n \big(\tilde\F^{^{\frak X}}\!, \Ker (\bar\rho^{^{\Omega,\frak X}})\big)\cong 
\Bbb H^n \big(\tilde\F^{^{\frak X,U}}\!, \Ker (\bar\rho^{^{\Omega,\frak X,U}})\big)
\eqno £5.4.4.$$

\medskip
£5.5. Let us  explicit the announced $\Ker (\bar\rho^{^{\Omega,\frak X}})\-$valued $2\-$cocycle.
For any $\F^{^{\frak X}}\!\-$morphism $\varphi\,\colon R\to Q\,,$ choose a lifting $x_\varphi$ 
in~$\M^{^{\Omega,\frak X}}\! (Q,R)$ (cf.~£4.3) and denote by $\bar x_\varphi$ the image of 
$x_\varphi$ in $\bar\M^{^{\Omega,\frak X}}\! (Q,R)\,;$ actually, we can do our choice in such a way that we have (cf.~£4.3)
$$\bar x_{\kappa^{_{\frak X}}_{_Q}(u)\circ\varphi} = 
\bar\upsilon^{_{\Omega,\frak X}}_{_Q} (u)\.\bar x_\varphi
\eqno £5.5.1\phantom{.}$$
for any $u\in Q\,,$ where $\kappa^{_{\frak X}}_{_Q}(u)\in \F^{^{\frak X}} \!(Q)$ denotes the conjugation by the image of~$u\,;$ indeed, if we have $\kappa^{_{\frak X}}_{_Q}(u)\circ\varphi =Ê\varphi$ then we get $u = \varphi (z)$ for a suitable $z\in Z(R)\,;$ since 
$\bar\M^{^{\Omega,\frak X}}$~is~{\it coherent\/}, in this case we obtain
$$\bar\upsilon^{_{\Omega,\frak X}}_{_Q} (u)\.\bar x_\varphi = 
\bar\upsilon^{_{\Omega,\frak X}}_{_Q} \big(\varphi (z)\big)\.\bar x_\varphi
= \bar  x_\varphi\.\bar\upsilon^{_{\Omega,\frak X}}_{_R} (z) = \bar x_\varphi
\eqno £5.5.2.$$
More precisely, if $Q$ and $R$ are contained in $P^{^U}$ and $\varphi\,\colon R\to Q$ comes from an $\F^{^{\frak X,U}}\!\-$morphism, it is quite clear that we may assume that $x_\varphi$ belongs to
$\big(N_{\! \M^{^{\Omega,\frak X}}} (U)\big)(Q,R)$ and then that 
$\frak h^{^{\frak X,U}}_{_\Omega} (x_\varphi)$  belongs to the image 
of~$\T_{\! L^{^U}_\F}^{^{\frak X^{^U}}} \!(Q,R)$ {\it via\/} $\frak l_{_{\F}}^{^{\frak X,U}}\,,$ so that actually we have (cf.~£5.3.6)
$$\bar\frak h^{^{\frak X,U}}_{_\Omega} (\bar x_\varphi) = \bar\frak m_{_{\F}}^{^{\frak X,U}} (\varphi)
\eqno £5.5.3.$$

\medskip
£5.6. Then, for any triple of subgroups $Q\,,$ $R$ and $T$ in~$\frak X\,,$ 
and any pair of $\F\-$morphisms $\psi\,\colon T\to R$ and $\varphi\,\colon R\to Q\,,$ since 
$ x_\varphi\. x_\psi$ and $ x_{\varphi\circ\psi}$ have the same image $\varphi\circ\psi$ in $\F(Q,T)\,,$ 
the {\it divisibility\/} of $\M^{^{\Omega,\frak X}}$ guarantees the existence and the uniqueness of
$k_{\varphi,\psi}\in {\rm Ker}(\rho^{_{\Omega,\frak X}}_T)$ fulfilling
$$ x_\varphi\. x_\psi =  x_{\varphi\circ\psi}\. k_{\varphi,\psi}
\eqno £5.6.1.$$
Denote by $\bar k_{\varphi,\psi}$ the image of $k_{\varphi,\psi}$ in 
${\rm Ker}(\bar\rho^{_{\Omega,\frak X}}_T)\,;$ since $\bar\M^{^{\Omega,\frak X}}\!$ is  
{\it coherent\/}, on the one hand for any $u\in Q$ and any $v\in R$ we get (cf.~£5.5.1)
$$\eqalign{\bar x_{\kappa^{_{\frak X}}_{_Q}(u)\circ\varphi}\. 
\bar x_{\kappa^{_{\frak X}}_{_R}(v)\circ\psi}
&= \big(\bar\upsilon^{_{\Omega,\frak X}}_{_Q} (u)\. \bar x_\varphi\big)\. 
\big(\bar\upsilon^{_{\Omega,\frak X}}_{_R} (v)\. \bar x_\psi\big)\cr
&= \bar\upsilon^{_{\Omega,\frak X}}_{_Q} \big(u\varphi (v)\big)\.\bar x_\varphi\.\bar x_\psi\cr
\bar  x_{(\kappa^{_{\frak X}}_{_Q}(u)\circ\varphi)\circ(\kappa^{_{\frak X}}_{_R}(v)\psi)} 
&= \bar x_{\kappa^{_{\frak X}}_{_Q} (u\varphi (v))\circ\varphi\circ\psi} = 
\bar\upsilon^{_{\Omega,\frak X}}_{_Q} \big(u\varphi (v)\big)\. \bar x_{\varphi\circ\psi}\cr}
\eqno £5.6.2;$$
hence, from the {\it divisibility\/} of $\bar\M^{^{\Omega,\frak X}}$ we obtain 
$$\bar k_{\kappa^{_{\frak X}}_{_Q}(u)\circ\varphi,
\kappa^{_{\frak X}}_{_R}(v)\circ\psi} = \bar k_{\varphi,\psi}
\eqno £5.6.3.$$
\eject
\noindent
That is to say,  for any $n\in \Bbb N\,,$ setting [3,~1.5]
$$\Bbb C^n \big(\tilde\F^{^{\Omega,\frak X}}\!,\Ker (\bar\rho^{_{\Omega,\frak X}})\big) = 
\prod_{\tilde\frak q\in \Fct(\Delta_n,\tilde\F^{^{\frak X}})}
{\rm Ker}(\bar\rho^{^{\Omega,\frak X}}_{\tilde\frak q (0)})
\eqno £5.6.4,$$
 we have obtained an element 
 $\bar k = \{\bar k_{\tilde\frak q}\}_{\tilde\frak q\in \Fct(\Delta_2,\tilde\F^{^{\frak X}})}$ in 
 $\Bbb C^2 \big(\tilde\F^{^{\Omega,\frak X}}\!,\Ker (\bar\rho^{_{\Omega,\frak X}})\big)$ where we set 
$\bar k_{\tilde\frak q} =\bar  k_{\tilde\frak q(1\bullet 2),\tilde\frak q (0\bullet 1)} = \bar  k_{\frak q(1\bullet 2),\frak q (0\bullet 1)}$ for some representative 
$\frak q\,\colon \Delta_2\to \F^{^{\frak X}}$ of $\tilde\frak q\,.$

\medskip
£5.7. We claim that $\bar k$ is actually a {\it $2\-$cocycle\/}; explicitly, considering the usual differential map [3,~A13.11]
$$\bar d^{^{\Omega,\frak X}}_2 : \Bbb C^2 \big(\tilde\F^{^{\frak X}}\!\!,
\Ker (\bar\rho^{^{\Omega,\frak X}})\big)\too 
\Bbb C^3 \big(\tilde\F^{^{\frak X}}\!\!,\Ker (\bar\rho^{^{\Omega,\frak X}})\big)
\eqno £5.7.1,$$
we claim that $\bar d^{^{\Omega,\frak X}}_2 (\bar k) = 0\,;$ indeed, with the notation above, for a third 
$\F^{^{\frak X}}\-$morphism $\eta\,\colon W\to T$ we get
$$\eqalign{(\bar x_\varphi\.\bar x_\psi)\.\bar x_\eta 
&= (\bar x_{\varphi\circ\psi} \.\bar k_{\tilde\varphi,\tilde\psi})\.\bar x_\eta
= (\bar x_{\varphi\circ\psi}\.\bar x_\eta)\.\big(\Ker (\bar\rho^{^{\Omega,\frak X}})
(\tilde\eta)\big) (\bar k_{\tilde\varphi,\tilde\psi})\cr
&= \bar x_{\varphi\circ\psi\circ\eta}\.\bar k_{\tilde\varphi\circ\tilde\psi,\tilde\eta}
\.\big(\Ker (\bar\rho^{^{\Omega,\frak X}})(\tilde\eta)\big)
(\bar k_{\tilde\varphi,\tilde\psi})\cr
\bar x_\varphi\.(\bar x_\psi\.\bar x_\eta) &= \bar x_\varphi\.(\bar x_{\psi\circ\eta}\.\bar k_{\tilde\psi,\tilde\eta}) = 
\bar x_{\varphi\circ\psi\circ\eta}\. \bar k_{\tilde\varphi,\tilde\psi\circ\tilde\eta}\.
\bar k_{\tilde\psi,\tilde\eta}\cr}
\eqno £5.7.2\phantom{.}$$
and   the {\it divisibility\/} of $\bar\M^{^{\frak X}}$ forces
$$\bar k_{\tilde\varphi\circ\tilde\psi,\tilde\eta}
\.\big(\Ker (\bar\rho^{^{\Omega,\frak X}})(\tilde\eta)\big)
(\bar k_{\tilde\varphi,\tilde\psi}) = \bar  k_{\tilde\varphi,\tilde\psi\circ\tilde\eta}
\. \bar k_{\tilde\psi,\tilde\eta}
\eqno £5.7.3;$$
since ${\rm Ker}(\bar\rho^{^{\Omega,\frak X}})$ is Abelian, in the additive notation we obtain
$$0 = \big(\Ker (\bar\rho^{^{\Omega,\frak X}})(\tilde\eta)\big)(\bar k_{\tilde\varphi,\tilde\psi}) -  \bar k_{\tilde\varphi,\tilde\psi\circ\tilde\eta} 
+ \bar k_{\tilde\varphi\circ\tilde\psi,\tilde\eta} - \bar  k_{\tilde\psi,\tilde\eta}
\eqno £5.7.4,$$
proving our claim.

\medskip
£5.8. Then,  in order to prove the existence of a section $\bar\sigma^{_{\Omega,\frak X}}\,,$
it suffices to show that $\bar k$ is a {\it $2\-$coboundary\/} and therefore, according to 
isomorphism~£5.4.4 above, it suffices to prove that the image {\it via\/} 
$\nu_{\bar\frak h^{^{\frak X,U}}_{_\Omega}}\!$ (cf.~£5.4.2) of the restriction of $\bar k$ to~$\tilde\F^{^{\frak X,U}}$ is a 
{\it $2\-$coboundary\/}. But, for any pair of  $\F^{^{\frak X,U}}\!\-$morphisms 
$\varphi\,\colon R\to Q$ and $\psi\,\colon T\to R\,,$ we have chosen $x_\varphi$ in 
$\big(N_{\!\M^{^{\Omega,\frak X}}\!} (U)\big)(Q,R)\,,$  $x_\psi$~in 
 $\big(N_{\!\M^{^{\Omega,\frak X}}\!} (U)\big)(R,T)$ and $x_{\varphi\circ\psi}$ in 
 $\big(N_{\!\M^{^{\Omega,\frak X}}\!} (U)\big)(Q,T)\,,$ so that in equality~£5.6.1 the element 
$k_{\varphi,\psi}$ belongs to  $\big(N_{\!\M^{^{\Omega,\frak X}}\!} (U)\big)(T)$ and therefore 
  we get (cf.~£5.3.4)
 $$\frak h^{^{\frak X,U}}_{_\Omega}  (x_\varphi)\.\frak h^{^{\frak X,U}}_{_\Omega} (x_\psi) =
 \frak h^{^{\frak X,U}}_{_\Omega}  (x_{\varphi\circ \psi})\. \frak h^{^{\frak X,U}}_{_\Omega} 
 (k_{\varphi,\psi})
 \eqno £5.8.1\phantom{.}$$
  and therefore we still get
 $$\bar\frak h^{^{\frak X,U}}_{_\Omega} (\bar x_\varphi)\.
 \bar\frak h^{^{\frak X,U}}_{_\Omega}(\bar x_\psi) = 
 \bar\frak h^{^{\frak X,U}}_{_\Omega} (\bar x_{\varphi\circ \psi})\. 
 \bar\frak h^{^{\frak X,U}}_{_\Omega}(\bar k_{\varphi,\psi})
 \eqno £5.8.2,$$
 \eject
 \noindent
 so that equalities~£5.5.3 force $\bar\frak h^{^{\frak X,U}}_{_\Omega} (\bar k_{\varphi,\psi}) = 1\,;$
that is to say,   the image {\it via\/} $\nu_{\bar\frak h^{^{\frak X,U}}_{_\Omega}}$ of the restriction of $\bar k$ to~$\tilde\F^{^{\frak X,U}}$ is just trivial, proving that $\bar k$ is a {\it $2\-$coboundary\/}.

\medskip
£5.9. Thus, we have obtained a {\it functorial section\/} $\bar\sigma^{_{\Omega,\frak X}}\,\colon 
\F^{^{\frak X}}\to \bar\M^{^{\Omega,\frak X}}$ of~$\bar\rho^{^{\Omega,\frak X}}\,;$ actually, $\bar\sigma^{_{\Omega,\frak X}}$ can be  modified  in order to get an {\it $\F^{^\frak X}\!\-$locality functorial section\/} [5,~2.9]. Indeed, for any $\F^{^\frak X}_{\! P}\-$morphism 
$\zeta\,\colon R\to Q\,,$ choosing $u_\zeta$ in $T_P (R,Q)$ lifting $\zeta\,,$
both $\bar\M^{^{\Omega,\frak X}}\-$morphisms $\bar\sigma^{_{\Omega,\frak X}}_{_{Q,R}} (\zeta)$ and $\bar\upsilon^{^{\Omega,\frak X}}_{_{Q,R}}(u_\zeta)$ (cf.~£4.3)
 lift $\zeta\,;$ once again,  the {\it divisibility\/} of~$\bar\M^{^{\Omega,\frak X}}$ guarantees the existence and the uniqueness  of~$\bar m_\zeta\in {\rm Ker}(\bar\rho^{^{\Omega,\frak X}}_{_R})$ fulfilling
$$\bar\upsilon^{^{\Omega,\frak X}}_{_{Q,R}}(u_\zeta) = 
\bar\sigma^{_{\Omega,\frak X}}_{_{Q,R}} (\zeta)\.\bar m_\zeta
\eqno £5.9.1;$$
actually, it follows easily from~£5.5.1 that $\bar m_\zeta$ only depends on 
$\tilde\zeta\in \tilde\F_{\!P}(Q,R)$ and, as above, we write $\bar m_{\tilde\xi}$ instead of $m_\xi\,;$ moreover, for a second $\F^{^\frak X}_{\! P}\-$mor-phism $\xi\,\colon T\to R\,,$ we get
$$\eqalign{\bar\sigma^{_{\Omega,\frak X}}_{_{Q,T}} (\zeta\circ\xi)\.\bar m_{\tilde\zeta\circ\tilde\xi} 
&= \bar\upsilon^{^{\Omega,\frak X}}_{_{Q,T}}(u_{\zeta\circ\xi}) = 
\bar\upsilon^{^{\Omega,\frak X}}_{_{Q,R}} (u_\zeta)\.
\bar\upsilon^{^{\Omega,\frak X}}_{_{R,T}}(u_\xi)\cr
&= \bar\sigma^{_{\Omega,\frak X}}_{_{Q,R}} (\zeta)\.\bar m_{\tilde\zeta}\.
\bar\sigma^{_{\Omega,\frak X}}_{_{R,T}} (\xi)\.\bar m_{\tilde\xi}\cr
&= \bar\sigma^{_{\Omega,\frak X}}_{_{Q,T}} (\zeta\circ\xi)\.\big(\Ker (\bar\rho^{^{\Omega,\frak X}})
(\tilde\xi)\big)(\bar m_{\tilde\zeta})\.\bar m_{\tilde\xi}\cr}
\eqno £5.9.2.$$

\smallskip
£5.10. Then, always the {\it divisibility\/} of $\bar\M^{^{\Omega,\frak X}}$  forces
$$\bar m_{\tilde\zeta\circ\tilde\xi} = \big(\Ker (\bar\rho^{^{\frak X}})
(\tilde\xi)\big)(\bar m_{\tilde\zeta})\.\bar m_{\tilde\xi}
\eqno £5.10.1\phantom{.}$$
and, since ${\rm Ker}(\bar\rho^{^{\Omega,\frak X}}_T)$ is Abelian (cf.~Proposition~£3.4,~£3.5 and~£3.6), in the additive notation we obtain
$$0 = \big(\Ker (\bar\rho^{^{\Omega,\frak X}})(\tilde\xi)\big)(\bar m_{\tilde\zeta}) - 
\bar m_{\tilde\zeta\circ\tilde\xi}+ \bar m_{\tilde\xi}
\eqno £5.10.2;$$
that is to say,  denoting by $\tilde\iota^{_\frak X}_P\,\colon \tilde\F^{^\frak X}_{\!P}\i 
\tilde\F^{^\frak X}$ the obvious inclusion functor, the correspondence $\bar m$ sending any  
$\tilde\F^{^{\frak X}}_{\! P}\-$morphism $\tilde\zeta\,\colon R\to Q$ 
to~$\bar m_{\tilde\zeta}$ defines a {\it $1\-$cocycle\/} in $\Bbb C^1
\big(\tilde\F^{^{\frak X}}_{\! P},\Ker (\bar\rho^{^{\Omega,\frak X}}) \circ\tilde\iota^{_\frak X}_P\big)\,;$ but, since the category $\tilde\F^{^\frak X}_{\! P}$ obviously
has a final object, we actually have [3,~Corollary~A4.8]
$$\Bbb H^1\big(\tilde\F^{^{\frak X}}_{\! P},\Ker (\bar\rho^{^{\Omega,\frak X}})
\circ\tilde\iota^{_\frak X}_P\big) = \{0\}
\eqno £5.10.3;$$
consequently, we obtain $\bar m = d^{^{\Omega,\frak X}}_0 (\bar w)$ for some element 
$\bar w = (\bar w_Q)_{Q\in \frak X}$ in 
$$\Bbb C^0\big(\tilde\F^{^\frak X}_{\! P},\Ker (\bar\rho^{^{\Omega,\frak X}})
\circ\tilde\iota^{_\frak X}_P\big) = \Bbb C^0\big(\tilde\F^{^\frak X},
\Ker (\bar\rho^{^{\Omega,\frak X}})\big)
\eqno £5.10.4.$$
In conclusion, equality~£5.9.1 becomes
$$\bar\upsilon^{^{\Omega,\frak X}}_{_{Q,R}}(u_\zeta) = 
\bar\sigma^{_{\Omega,\frak X}}_{_{Q,R}} (\zeta)\.\big(\Ker (\bar\rho^{^{\Omega,\frak X}})
(\tilde\zeta)\big)(\bar w_Q)\.\bar w_R^{-1} = \bar w_Q\.\bar\sigma^{_{\Omega,\frak X}}_{_{Q,R}} (\zeta)\.\bar w_R^{-1}
\eqno £5.10.5;$$
\eject
\noindent
thus,  the new correspondence which, for any pair of subgroups $Q$ and $R$ in $\frak X\,,$  sends any $\varphi\in \F (Q,R)$ to  $\bar w_Q\.\bar\sigma^{_{\Omega,\frak X}}_{_{Q,R}} 
(\varphi)\.\bar w_R^{-1}$  defines  an {\it $\F^{^{\frak X}}\!\-$locality functorial section\/}   
of $\bar\rho^{^{\frak X}}\,.$ We are done.

\bigskip
\centerline{\large References}
 
\bigskip

\smallskip\noindent
[1] Bob Oliver, {\it A remark on the construction of centric linking systems\/}, 
arxiv.org/abs/1612.02132
\smallskip\noindent
[2]\phantom{.} Llu\'\i s Puig, {\it Frobenius categories},
Journal of Algebra, 303(2006), 309-357.
\smallskip\noindent
[3]\phantom{.} Llu\'\i s Puig, {\it ``Frobenius categories versus Brauer blocks''\/}, Progress in Math. 
274(2009), Birkh\"auser, Basel.
\smallskip\noindent
[4]\phantom{.} Llu\'\i s Puig, {\it A criterion on trivial homotopy\/},  arxiv.org/abs/1308.3765.
\smallskip\noindent
[5]\phantom{.} Llu\'\i s Puig, {\it Existence, uniqueness and functoriality of the perfect locality over a Frobenius  $P\-$category\/}, Algebra Colloquium, 23(2016) 541-622.

\end